\documentclass[11pt,twoside,a4paper]{article}
\pdfoutput=1
\usepackage[T1]{fontenc}
\usepackage{hyperref}
\usepackage{amsmath,amsfonts,amssymb,amscd,amsthm}
\usepackage[latin1]{inputenc}
\usepackage[numbers]{natbib}
\usepackage{color}
\usepackage[usenames,dvipsnames]{xcolor}
\usepackage{doi}
\usepackage{graphicx}
\usepackage{booktabs}
\usepackage[lighttt]{lmodern}
\usepackage[labelfont={bf,small},font=small]{caption}
\usepackage[vlined]{algorithm2e}
\SetKwProg{Fn}{Procedure}{}{end}
\SetAlCapSkip{0.5em}

\SetAlgoCaptionLayout{SetSmall}

\hypersetup{
    bookmarks=true,
    breaklinks=true,
    bookmarksopen=true,
    pdftitle={Unified framework for the efficient solution of n-field coupled problems},    
    pdfauthor={F. Verdugo and W. A. Wall},     
    colorlinks=true,       
    linkcolor=black,          
    citecolor=black,        
    filecolor=black,      
    urlcolor=black           
}

\usepackage{titlesec}

\usepackage{mydefs}

\usepackage{fancyhdr}
\colorlet{rev1}{black}
\colorlet{rev2}{black}

\newcommand{\Partial}[2]{\dfrac{\partial #1}{\partial #2}}

\begin{document}

\thispagestyle{empty}

{ \LARGE \bf Unified computational framework for the efficient solution\\ of n-field coupled problems with monolithic schemes\footnotemark}

\footnotetext{Article accepted for publication in \emph{Computer Methods in Applied Mechanics and Engineering}, July 14 2016.}


Francesc Verdugo and Wolfgang A. Wall

\vspace{0.4cm}

Institute for Computational Mechanics\\
Technische Universit\"at M\"unchen\\
Boltzmannstr. 15\\ 
85748 Garching b.~M\"unchen\\
{Email: \texttt{\{verdugo,wall\}@lnm.mw.tum.de }}

%
%
%
%
%
%
%
%
%
%

\rule{\textwidth}{1pt}

\textbf{Abstract}

{\color{rev1} 
In this paper, we propose and evaluate the performance of a unified computational framework for preconditioning  systems of linear equations resulting from the solution of coupled problems with monolithic schemes. The framework is composed by promising application-specific preconditioners presented previously in the literature with the common feature that they are able to be implemented for a generic coupled problem, involving an arbitrary number of fields, and to be used to solve a variety of applications.
}
The first selected preconditioner is based on a generic block Gauss-Seidel iteration for uncoupling the fields, and standard algebraic multigrid (AMG) methods for solving the resulting uncoupled problems. The second preconditioner is based on the {\color{rev1} semi-implicit method for pressure-linked equations (SIMPLE)} which is extended here to deal with an arbitrary number of fields, {\color{rev1} and} also results  in uncoupled problems that can be solved with standard AMG. Finally, a more sophisticated preconditioner is considered which enforces the coupling at all AMG levels, in contrast to the other two techniques which resolve the coupling only at the finest level. {\color{rev1} Our purpose is to show that these methods perform satisfactory in quite different scenarios apart from their original applications. To this end, we consider}   three very different coupled problems:  thermo-structure interaction, fluid-structure interaction and a complex  model of the human lung. Numerical results show that {\color{rev1} these} general purpose methods are efficient and scalable in {\color{rev1} this  range of}  applications.

\textbf{Keywords:} Coupled problems $\cdot$ Monolithic solution schemes   $\cdot$   Block preconditioners  $\cdot$ Algebraic Multigrid  $\cdot$ High Performance Computing

\vspace{-0.3cm}
\rule{\textwidth}{1pt}

\section{Introduction}

The numerical simulation of coupled problems requires special solution strategies for solving the coupling between the underlying physical fields.
The so-called \emph{partitioned} methods  \citep{FelippaParkFarhat2001} are often the preferred choice in industrial applications  because they allow to reuse existing (black-box) solvers for the individual fields. However, this approach is unstable for many challenging strongly coupled problems 
  \citep{CausinGerbeauNobile2005,FoersterWallRamm2007} and, therefore, another family of methods called \emph{monolithic} are required in a variety of complex settings. It has also been shown that monolithic schemes are often preferable in terms of efficiency as compared to partitioned ones.
For that reason, monolithic methods have been the preferred option for solving many strongly coupled problems, see e.g. \citep{BadiaQuainiQuarteroni2008,BadiaMartinPlanas2014,CyrShadidTuminaroPawlowskiChacon2013,DanowskiGravemeierYoshiharaWall2013,GeeKuettlerWall2011,HeilHazelBoyle2008,HowleKirby2012,
HowleKirbyDillon2013,KuettlerGeeFoersterComerfordWall2010,LinShadidTuminaroSalaHenniganPawlowski2010}.

The price to be paid for the extra robustness of monolithic methods is a more challenging system of linear equations. The system matrix is a big sparse matrix with a special block structure representing each of the underlying physical fields and frequently has a very bad condition number. 
In real-world applications, iterative methods as GMRES \citep{SaadSchultz1986} are required to attack this linear system, which requires efficient preconditioners for addressing the bad conditioning of the problem. 
Selecting a suitable preconditioner is the key point in the solution process and it is crucial to benefit from the robustness  and efficiency of a monolithic approach.

The preconditioners proposed in the literature for coupled problems are usually complex methods designed for specific applications. This approach has the advantage that the preconditioner can be optimized for the particular problem leading to efficient solvers, but on the other hand, specific implementations have to be developed and maintained by experienced users for each particular application. This is specially an issue in ``general purpose'' multiphysics codes designed to simulate  different kinds of problems,  because specific implementations have to be developed and  maintained separately.

In this paper, we {\color{rev1} explore} an alternative approach. {\color{rev1} We select from the preconditioners available in the literature several techniques that are able to be implemented for a generic coupled problem, involving an arbitrary number of fields, and to be used to solve a variety of applications. We have implemented these methods in a unified computational framework allowing researchers to test and use them conveniently for different {\color{rev1} problem types} without the need of developing application-specific implementations. }Our purpose is to show that such a {\color{rev1} framework} is realizable and that {\color{rev1} the selected methods are} efficient and able to cope with a wide range of applications.

The first type of preconditioners proposed here is based on the so-called \emph{block preconditioners}, see e.g. \citep{HowleKirbyDillon2013}, which lead to efficient solvers by taking advantage of the block structure of the monolithic system. The key idea is to use approximate block inverses in order to untangle the coupling between the physical fields, and then, to use efficient algebraic multigrid (AMG) \citep{BriggsHensonMcCormick2000,TrottenbergOosterleeSchueller2001,VanekMandelBrezina1996} solvers for the resulting uncoupled problems. This approach has been widely used in the literature for particular applications as fluid-structure interaction (FSI) \citep{GeeKuettlerWall2011,KuettlerGeeFoersterComerfordWall2010}, magneto hydro dynamics  \citep{BadiaMartinPlanas2014,CyrShadidTuminaroPawlowskiChacon2013,ElmanHowleShadidShuttleworthTuminaro2008} and thermo-structure interaction (TSI) \citep{DanowskiGravemeierYoshiharaWall2013}, among others. In this context, we propose two different methods: a generalization of our previous works \citep{DanowskiGravemeierYoshiharaWall2013,GeeKuettlerWall2011}
considering a block Gauss-Seidel (BGS) iteration for uncoupling the fields, and a generalization of \citep{BadiaMartinPlanas2014} which considers an extension of the SIMPLE method \citep{CarettoGosmanPatankarSpalding1973} {\color{rev1} for the same purpose}.

The drawback of standard block preconditioning {\color{rev1} based on BGS or SIMPLE} is that the coupling is resolved only  at the finest multigrid level. Thus, using efficient AMG solvers for the {\color{rev1} underlying problems} does not necessarily imply a good treatment of the coupling and a fast global solution.  This drawback is overcome by \citet{GeeKuettlerWall2011} who propose an enhanced block preconditioner (referred to as \emph{monolithic AMG}) for FSI applications, which enforces the coupling at all multigrid levels. {\color{rev1} This often  results} in  a better performance of the solver {\color{rev1} in the FSI problems studied in \citep{GeeKuettlerWall2011}}. 
Monolithic AMG preconditioners are a promising approach for other coupled problems as well but,  to our knowledge, this strategy has only been applied to FSI so far.  One of the reasons why this approach has not been used and transferred to other problems so far is exactly due to the reasons given above: the required special know-how and cumbersome, necessary separate implementation. Here, we also propose an extension of the monolithic AMG preconditioner to a generic coupled problem. This is the first time that this preconditioner type is transferred to other problem types than the original FSI-application. 
{\color{rev2} Other more sophisticated block preconditioners (not discussed here) can also circumvent the degradation in performance associated with standard block preconditioners based on BGS or SIMPLE. For instance, the so-called pressure-convection-diffusion (PCD) or the least square commutator (LSC) preconditioners are successful alternatives in the context of the Navier-Stokes equations \citep{CyrShadidTuminaro2012}. For  magneto hydro dynamics (MHD), suitable alternatives are the preconditioners based on operator splitting or approximate commutators as shown in \citep{CyrShadidTuminaroPawlowskiChacon2013,PhillipsElmanCyrShadidPawlowski2014} respectively.}

{\color{rev2} The block preconditioners discussed in this paper} {\color{rev1} have been implemented in the high performance multiphysics code BACI \citep{WallGeeBaci} and have been tested with three very different applications:  TSI, FSI and a complex respiratory model of the human lung proposed by \citet{YoshiharaIsmailWall2013,YoshiharaRothWall2015}. The numerical results show that these type of preconditioners  show good performance for these different applications and suggest that the methods can be potentially applied to other coupled problems as well.
}

{\color{rev2} The block preconditioners studied here have been also implemented up to a certain extent in other related software packages and successfully applied to model problems different from those considered in this paper. One of the first object-oriented libraries supporting block preconditioners was already presented back in 1998 \citep{ChowHeroux1998}. This package included matrices with block structure and classical methods for uncoupling the fields as block Jacobi or block Gauss-Seidel. However, it did not support more advanced features as recursively nested block matrices which where not introduced until recently in modern software tools as Playa \citep{HowleKirbyLongBrennanKennedy2012}, Teko \citep{CyrShadidTuminaro2016,TekoWeb} or PCFIELDSPLIT \citep{BrownKnepleyMayMcInnesSmith2012,pcfieldsplitWeb}. Playa is a general high-level interface for manipulating block matrices which in particular allows the creation of block preconditioners. On the other hand, Teko is a block preconditioning environment that includes state-of-the-art preconditioners for Navier-Stokes \citep{CyrShadidTuminaro2012} and MHD  \citep{PhillipsElmanCyrShadidPawlowski2014} ready to use, modify and extend by the users. Both Teko and Playa are part of the Trilinos project \cite{TrilinosWeb}.  Finally, PCFIELDSPLIT is a component of the PETSc library \citep{PETSCWeb} and shares some similarities with Teko like the creation and run-time design of complex block preconditioners including modern methods based on Schur complements. 

The computational framework discussed in this paper shares many of the capabilities of Teko and PCFIELDSPLIT including the support for nested block matrices and run-time design of preconditioners. In contrast to these two packages, we have included and tested in our framework the  advanced monolithic AMG solver where a multilevel hierarchy of block matrices has to be created. Moreover, we have provided the ability of combining monolithic AMG solvers
 with incomplete block factorizations based on the Schur complement as it is shown in Section 6 for a pulmonary mechanics example.}


The remainder of the paper is structured as follows. In section 2, the monolithic approach is presented for the solution of a generic coupled problem leading to the system of linear equations to be solved with the proposed techniques. Section 3 includes a brief overview of solution techniques for linear systems and describes the role of the preconditioner. In section 4, the proposed general framework for preconditioning coupled problems is presented. In Section 5, some remarks on the implementation  are given.
Section 6 details the three applications considered for demonstrating the versatility of the methods {\color{rev2} and Section 7 contains the corresponding numerical examples. Finally, the paper is concluded in Section 8.} 

\section{Monolithic system of linear equations}

We focus on the solution of coupled problems that are described by partial differential equations {\color{rev1} that} after the usual space and time discretization result in a set of non-linear algebraic equations to be solved at each time step. For this purpose, we consider a generic set of non-linear equations
\begin{equation}
\label{eq:non-linear-problems}
\begin{aligned}
\fmat_1(\umat_1,\umat_2,\ldots,\umat_N) &= \zeromat,\\
\fmat_2(\umat_1,\umat_2,\ldots,\umat_N) &= \zeromat,\\
&\vdots\\
\fmat_N(\umat_1,\umat_2,\ldots,\umat_N) &= \zeromat,
\end{aligned}
\end{equation}
representing the discrete equations associated with a generic coupled problem with $N$ coupled fields. The non-linear vector-valued functions  $\fmat_1,\ldots,\fmat_N$
 represent the discrete versions of the state equations for each of the  fields, and the vectors $\umat_1,\umat_2,\ldots,\umat_N$ stand for the {\color{rev1}corresponding} discrete solution vectors. Particular coupled problems  leading to systems like~\eqref{eq:non-linear-problems} are given in Sections \ref{sec:TSI-theory}, \ref{sec:Intro:FSI} and \ref{sec:lung:theory} for the applications considered in the paper.

The discrete equations \eqref{eq:non-linear-problems} form a fully-coupled non-linear problem which is solved here using a monolithic Newton method.  To facilitate the discussion,  equations \eqref{eq:non-linear-problems} are written more compactly as
\begin{equation}
\label{eq:compact-nlin-problem}
\Fmat(\Umat)=\zeromat
\quad\text{with}\quad
\Fmat:=
\left[
\begin{array}{c}
\fmat_1\\
\fmat_2\\
\vdots\\
\fmat_N
\end{array}
\right]
\quad\text{and}\quad
\Umat:=
\left[
\begin{array}{c}
\umat_1\\
\umat_2\\
\vdots\\
\umat_N
\end{array}
\right].
\end{equation}
With this notation, the monolithic Newton iteration for solving problem \eqref{eq:non-linear-problems} is:
\begin{equation}
\Umat^{(i+1)}=\Umat^{(i)} + \Delta\Umat^{(i)}
\end{equation}
where $\Umat^{(i)}$ stands for the full solution vector at the $i$-th iteration, and the correction $\Delta\Umat^{(i)}$ is obtained by solving the system of linear equations
\begin{equation}
\label{eq:compact-system}
\left[\dfrac{\partial\Fmat}{\partial\Umat}\right]^{(i)}\Delta\Umat^{(i)} = -\Fmat^{(i)}.
\end{equation}

System  \eqref{eq:compact-system} results from the usual linearization of the non-linear problem \eqref{eq:non-linear-problems} at the $i$-th Newton step and it involves 
the linearization of all the field       equations $\fmat_1,\fmat_2,\ldots,\fmat_N$ with respect to all the variables $\umat_1,\umat_2,\ldots,\umat_N$. This leads to a system matrix with an $N\times N$ block structure:
\begin{equation}
\label{eq:linearized-coupled-problem}
\def\arraystretch{2.3}
\left[
\begin{array}{cccc}
\dfrac{\partial\fmat_1}{\partial\umat_1} & \dfrac{\partial\fmat_1}{\partial\umat_2} & \cdots &\dfrac{\partial\fmat_1}{\partial\umat_N} \\
\dfrac{\partial\fmat_2}{\partial\umat_1} & \dfrac{\partial\fmat_2}{\partial\umat_2} & \cdots &\dfrac{\partial\fmat_2}{\partial\umat_N} \\
\vdots & \vdots & \ddots & \vdots\\
\dfrac{\partial\fmat_N}{\partial\umat_1} & \dfrac{\partial\fmat_N}{\partial\umat_2} & \cdots &\dfrac{\partial\fmat_N}{\partial\umat_N} \\
\end{array}
\right]^{(i)}
\left[
\begin{array}{c}
\Delta\umat_1\\
\Delta\umat_2\\
\vdots\\
\Delta\umat_N
\end{array}
\right]^{(i)}
=
-
\left[
\begin{array}{c}
\fmat_1\\
\fmat_2\\
\vdots\\
\fmat_N
\end{array}
\right]^{(i)}.
\end{equation}
It is implicitly assumed that some of the blocks in \eqref{eq:linearized-coupled-problem} can also be zero.

The system \eqref{eq:linearized-coupled-problem} is one of the main difficulties associated with the solution of coupled problems with monolithic schemes.  From the linear solver perspective, system  \eqref{eq:linearized-coupled-problem}  is particularly challenging because it involves all the physical fields. As a result, the system size is especially large and  the condition number is usually very bad due to the different {\color{rev1} scaling magnitudes associated with the distinct} nature of the underlying physics.

The goal of this paper is to propose solution strategies, i.e. preconditioners, for the efficient solution of system \eqref{eq:linearized-coupled-problem} for a wide range of coupled problems.

\section{Iterative solution techniques}

Before presenting the preconditioners for coupled problems, {\color{rev1} since this paper is targeting people dealing with modeling and simulation of coupled problems which are not necessarily solver experts}, we briefly describe the linear solver used for {\color{rev1} solving} the monolithic system \eqref{eq:linearized-coupled-problem}, and which role the preconditioners  play within the method. Moreover, we briefly comment on multigrid preconditioning techniques for single-field problems which are the main ingredient  for building the proposed preconditioners for coupled problems. {\color{rev1}Obviously, it is beyond the scope of this paper to provide an adequate in-depth discussion of these topics. The idea is more to provide some intuition on the subject and some references for further consideration by the reader.}

\subsection{Krylov subspace methods}

The resolution of the monolithic system \eqref{eq:linearized-coupled-problem} with direct methods is not affordable in many practical applications. The fill-in of the underlying sparse matrices leads to extreme memory requirements and computation times making direct methods inefficient. Moreover, direct methods are  complex to implement in parallel and have a rather poor scalability. For this reason, the so-called  Krylov subspace methods \cite{Vorst2003} are often preferred for solving large sparse systems like \eqref{eq:linearized-coupled-problem} in parallel computers.

Among the most known variants of Krylov subspace methods are the conjugate gradient method \cite{Ginsburg1963} for symmetric positive definite matrices, and the GMRES method \citep{SaadSchultz1986} which is able to deal with generic non-singular matrices.
In the following, we consider the GMRES method  implemented in the AztecOO package \cite{Heroux2007} for solving the coupled system~\eqref{eq:linearized-coupled-problem}.  We choose GMRES because is the most general Krylov subspace method able to deal with non-symmetric matrices, but the same discussion applies for other Krylov solvers. 

For a given system of linear equations 
 \begin{equation}
 \label{eq:generic-system}
\Amat\xmat=\bmat,
\end{equation}
the GMRES method furnishes a sequence of approximate solutions $\xmat^k$ that converge to the exact solution $\xmat$ starting from an initial guess $\xmat^0$.
The convergence {\color{rev1} rate} strongly depends on the spectral properties of the system matrix $\Amat$. {\color{rev1} In the easiest case, when matrix $\Amat$ is symmetric definite,} the convergence {\color{rev1} rate} strongly depends on its condition number {\color{rev1}(see e.g. \cite{TrefethenBau1997}):}
\begin{equation*}
\kappa_2(\Amat) := \dfrac{\sigma_\mathrm{max}(\Amat)}{\sigma_\mathrm{min}(\Amat)},
\end{equation*}
where $\sigma_\mathrm{max}(\Amat)$ and $\sigma_\mathrm{min}(\Amat)$ are the maximum and minimum eigenvalue of matrix $\Amat$ respectively. If the condition number is small, i.e. $\kappa_2(\Amat)\approx 1$, then, the convergence is expected to be fast. Whereas, for large condition numbers, the convergence is frequently very slow. 
{\color{rev1} The connection between the eigenvalues of $\Amat$ and the convergence behavior of the Krylov method becomes more involved when the matrix is not symmetric definite. The condition number is only a simplified convergence estimate in these situations. For a normal matrix $\Amat$, the convergence depends on how the eigenvalues are distributed on the complex plain. If the eigenvalues are clustered within a small region  away from the origin,  fast convergence of the method is expected, see e.g. \citep{Saad2003}. Unfortunately, the convergence behavior of Krylov methods is still not fully understood for general nonsingular matrices which are not symmetric definite nor normal, see e.g. \citep{LiesenTichy2004}. Consequently, the condition number can only be used as a simplified convergence estimate for non symmetric definite matrices.}

The condition number of the monolithic matrix in \eqref{eq:linearized-coupled-problem}  is usually very large {\color{rev1} due to several reasons including different scaling magnitudes associated with different underlying physics, strong hyperbolic and elliptic operators, and fine finite element meshes. Efficient preconditioners have to be used for improving the spectral properties of the system matrix in order use Krylov methods efficiently.  Without the use of a preconditioner}, Krylov  methods do not converge at all for many practical applications. However, if efficient preconditioners are used, {\color{rev1} Krylov methods} are among the most efficient techniques for solving this type of systems, see for instance \citep{LinShadidTuminaroSalaHenniganPawlowski2010}.

\subsection{The idea of preconditioning}

A preconditioner $\Pmat$ for a generic system \eqref{eq:generic-system} is a linear operator that tries to represent the action of the inverse of the system matrix $\Amat$, namely $\Pmat\approx\Amat^{-1}$. In practice, the preconditioner is a procedure which takes a given vector $\bar\bmat\in\mathbb{R}^r$ ($r$ is the number of rows or columns in $\Amat$) and returns a vector  $\bar\xmat\in\mathbb{R}^r$ being the approximate solution of a system of equations formed by the system matrix $\Amat$ and taking vector $\bar\bmat$ as the right hand side: 
\begin{equation}
\begin{array}{llll}
\Pmat:& \mathbb{R}^r& \longrightarrow & \mathbb{R}^r\\
      & \bar{\bmat} & \longmapsto & \bar\xmat \approx\Amat^{-1}\bar{\bmat}
\end{array}.
\end{equation}
The preconditioner can be seen as a cheap linear solver used to improve the convergence of a GMRES iteration. In fact, the preconditioner is frequently referred to as \emph{the solver}, whereas the outer GMRES iteration as \emph{the accelerator}. This is because the preconditioner usually requires more attention and problem-specific fine tuning than the GMRES iteration which is a generic procedure that can be considered as a black box.

For a given preconditioner $\Pmat$, the preconditioned version of the system \eqref{eq:generic-system} is
\begin{equation*}
\left(\Amat\Pmat\right)\left(\Pmat^{-1}\xmat\right)=\bmat.
\end{equation*}
Solving the preconditioned system requires first solving 
\begin{equation}
\label{eq:right-preconditioned-system}
\left(\Amat\Pmat\right)\ymat=\bmat
\end{equation}
for the auxiliary variable $\ymat:=\Pmat^{-1}\xmat$, and finally obtaining the unknown $\xmat$ with an application of the preconditioning operator, $\xmat=\Pmat\ymat.$

If the preconditioner $\Pmat$ is a good approximation of $\Amat^{-1}$, then the matrix $\Amat\Pmat$ is expected to be close to the identity and the condition number of the preconitioned system \eqref{eq:right-preconditioned-system} is expected to be $\kappa_2(\Amat\Pmat^{-1})\approx 1 $ {\color{rev1} which is usually associated with} a good convergence of the GMRES method. The extra cost to be paid for this improved conditioning is an application of the preconditioner at each GMRES iteration, see \citep{SaadSchultz1986}.

Our purpose is to present different {\color{rev1}preconditioners} (i.e. efficient solvers) for the monolithic system \eqref{eq:linearized-coupled-problem} to be used within a GMRES iteration for solving the system efficiently. The key idea of the proposed methods is to reuse existing preconditioning techniques that are known to be efficient for the underlying single-field problems. In this context,  one of the most successful preconditioners are multigrid methods.

\subsection{Multigrid methods}
\label{sec:multigrid}


Multigrid methods \cite{BriggsHensonMcCormick2000} are efficient solution techniques designed for systems of linear equations resulting from the discretization of partial differential equations.
The basic idea is to {\color{rev1} solve} the system with an inexpensive iterative method (e.g. few Gauss-Seidel iterations), and then to approximate the unresolved error  using a coarser version of the problem (that can be constructed by discretizing the underlying PDE using a coarser mesh). This idea can be applied recursively to an arbitrary number of coarser levels until the system matrix is small enough to compute an accurate coarse correction using a direct solver.

{\color{rev1} Defining} a multigrid method requires introducing a hierarchy of coarser versions of the original discretization.  In this context, one can distinguish two main approaches: either geometric multigrid methods \citep{WesselingOosterlee2001}, which explicitly consider coarser computational meshes, or algebraic multigrid methods (AMG) \cite{BriggsHensonMcCormick2000}, which achieve coarser approximation spaces  precluding the computation of new meshes by using information of the system matrix. The chosen coarsening procedure has little algorithmic relevance here. 
In the examples below, we use the so-called aggregation-based AMG methods \cite{VanekBrezinaMandel2001}  provided in the packages MueLu \citep{ProkopenkoHuWiesnerSiefertTuminaro2014} and ML \cite{GeeSiefertHuTuminaroSala2006} because they are able to deal with complex 3D geometries.

\begin{figure}[ht!]
\centering
\includegraphics[width=0.6\textwidth]{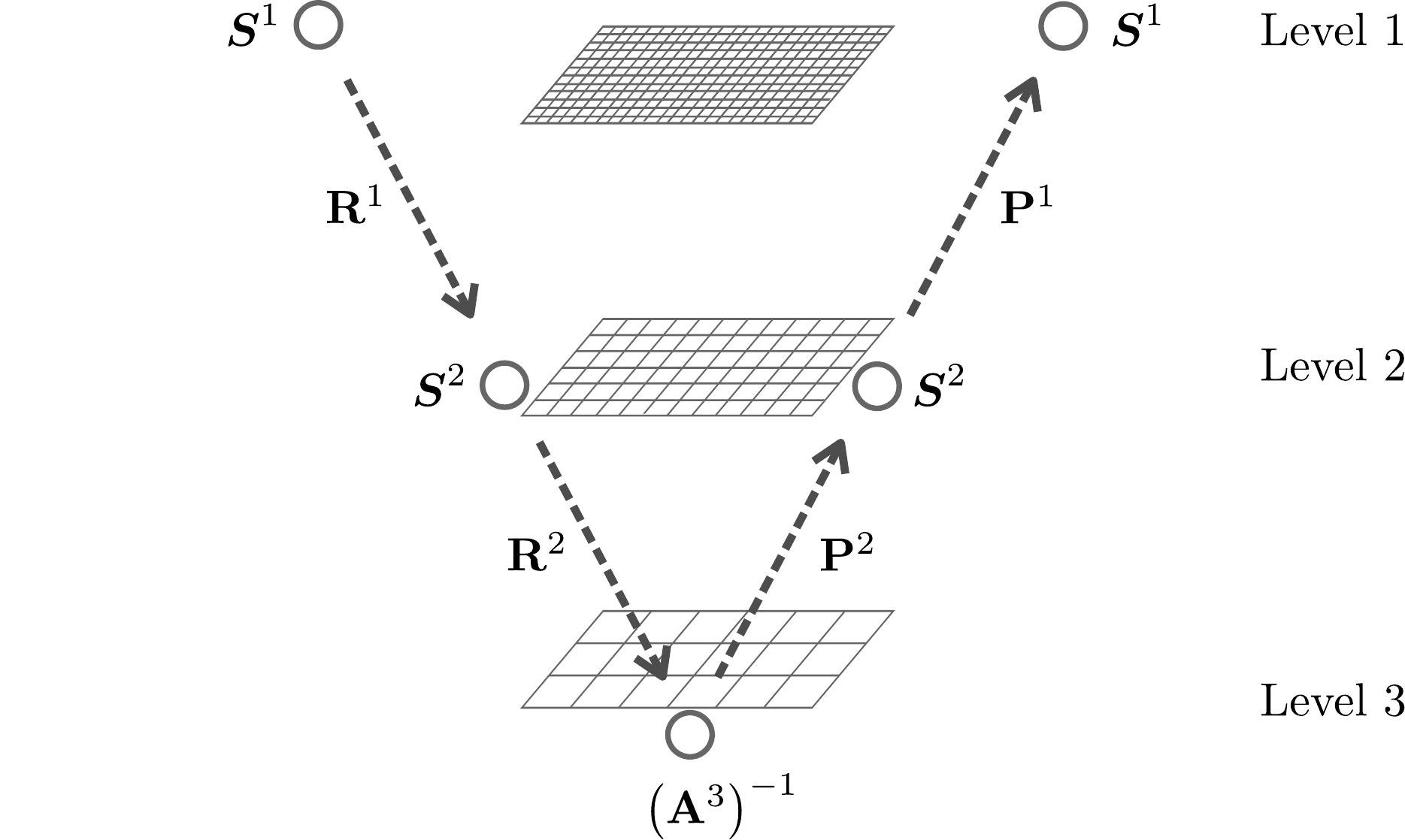}
\caption{Illustration of a multigrid V-cycle with $3$ levels.} 
\label{fig:V-cycle}
\end{figure}

The computed hierarchy of discretizations form the so-called multigrid levels. We denote the finest level as $\ell=1$, whereas  $\ell=\ell_\mathrm{max}$ corresponds to the coarsest level, see Figure \ref{fig:V-cycle}.
The main ingredients of a multigrid method are the following objects associated with each one of the levels  $\ell=1,\ldots,\ell_\mathrm{max}$:  1) the projection matrices $\Pmat^\ell$, that translate vectors form the coarse discretization on level $\ell+1$ to the fine discretization on level $\ell$, 2) Restriction matrices $\Rmat^\ell$,  that transforms vectors from level $\ell$ to level $\ell+1$, and 3) iterative methods $\Svec^\ell$ (referred to as \emph{the smoothers}) used to solve the system of equations  associated with the discretization on level $\ell$. 
Using the transfer operators, the matrices $\Amat^\ell$ associated with the discretization on levels $\ell=2,\ldots,\ell_\mathrm{max}$  are computed from the fine level matrix $\Amat^1$ as 
\begin{equation}
\Amat^{\ell+1}=\Rmat^\ell\Amat^\ell\Pmat^\ell,\quad\text{for }\ell=1,\ldots,\ell_\mathrm{max}-1.
\end{equation}
Note that with this formula, computing the coarse versions of matrix $\Amat^1$ do not require constructing coarser meshes and performing the standard assembly operations at each one of the levels. In other words, the multigrid method is a purely algebraic procedure, once the transfer operators $\Pmat^\ell$ and $\Rmat^\ell$ are available.

A multigrid method is presented in Algorithm \ref{alg:vcycle}.
The basic idea of the algorithm is to compute an approximation of the solution of the system of equations associated on the finest multigrid level $\ell=1$ using the smoother $\Svec^1$, i.e. a cheap iterative method that is very efficient in reducing high frequency error components, and  introducing an auxiliary system of equations on the next level $\ell=2$, where remaining error components are again highly oscillatory, for computing a correction. This auxiliary system on level $\ell=2$ is approximated following the same idea by using the smoother $\Svec^2$ and introducing another auxiliary system at the next level $\ell= 3$. This coarsening procedure is repeated until   the coarsest level $\ell=\ell_\mathrm{max}$ is reached. At the coarsest level, the problem can be efficiently {\color{rev1} solved} with a direct {\color{rev1} method}, and therefore, no further coarsening is needed. When the correction on the coarsest level is available, it is projected back through the different levels until the finest discretization is reached. The resulting method is called a V-cycle denoting the sequential order in which the different levels are visited, see  Figure \ref{fig:V-cycle}.

\begin{algorithm}
\caption{Multigrid V-cycle with $\ell_\mathrm{max}$ levels for the resolution of the fine level system $\Amat^1\xmat^1=\bmat^1$. }
\label{alg:vcycle}
\Fn{$\boldsymbol{\Phi}^\mathrm{V-cycle}(\Amat^\ell,\bmat^\ell,\xmat^\ell,\ell)$}{
\KwIn{System matrix $\Amat^\ell$, right hand side $\bmat^\ell$, initial guess $\xmat^\ell$, and current level $\ell$.}
\KwOut{Approximated solution of the system $\Amat^\ell\xmat^\ell=\bmat^\ell$.}
\eIf{$\ell<\ell_\mathrm{max}$ \textit{(fine level)}}{
Apply pre-smoother:\\
$\xmat^\ell \leftarrow \Svec^\ell(\Amat^\ell,\bmat^\ell,\xmat^\ell)$\\
{Compute coarse level matrix:}\\
$\Amat^{\ell+1}=\Rmat^\ell\Amat^\ell\Pmat^\ell$\\
{Compute coarse level residual:}\\
$\bmat^{\ell+1}=\Rmat^\ell(\bmat^\ell-\Amat^\ell\xmat^\ell)$\\
{Compute coarse level correction recursively:}\\
$\xmat^{\ell+1} = \zeromat$\\
$\xmat^{\ell+1}\leftarrow\boldsymbol{\Phi}^\mathrm{V-cycle}(\Amat^{\ell+1},\bmat^{\ell+1},\xmat^{\ell+1},\ell+1)$\\
{Apply correction:}\\
$\xmat^\ell \leftarrow \xmat^\ell + \Pmat^\ell\xmat^{\ell+1}$\\
{Apply post-smoother:}\\
$\xmat^\ell \leftarrow \Svec^\ell(\Amat^\ell,\bmat^\ell,\xmat^\ell)$\\
\KwRet $\xmat^\ell$\\
}(\textit{(coarsest level)})
{
{Apply a direct solver:}\\
$\xmat^\ell\leftarrow\left(\Amat^\ell\right)^{-1}\bmat^\ell$\\
\KwRet $\xmat^\ell$\\
}
 }
\end{algorithm}

Multigrid techniques are well established mainly for single-field problems. For elliptic problems (e.g. the {\color{rev1}steady}-state heat equation or structural mechanics),  multigrid methods are optimal in the sense that the cost of solving the system is proportional to the system size, cf. \citep{BrennerScott2008}, and can be implemented in parallel rendering good {\color{rev1} scalability}.
The application of standard multigrid techniques (i.e. aggregation-based methods \cite{VanekBrezinaMandel2001})  to coupled problems is more challenging. 
Using standard aggregation-based AMG in coupled problems is known  as \emph{fully coupled} AMG and it is efficient for a number of applications including  incompressible flows coupled with transport  \citep{LinSalaShadidTuminaro2006},  resistive magneto hydrodynamics \citep{LinShadidTuminaroSalaHenniganPawlowski2010}, {\color{rev1} and semiconductor device modeling \citep{LinShadidSalaEtAl2009}}. However, the applicability of fully coupled AMG as a general purpose preconditioner for coupled problems is limited. 
This approach requires that the degrees of freedom of the underlying physical fields are collocated at the same mesh nodes, which is not fulfilled in many cases. {\color{rev1} That is, the method requires the same finite element interpolation for all field variables.}
Thus, more versatile methods are required for preconditioning coupled problems.

\section{Unified framework for preconditing coupled problems}
\label{sec:all-precs}

In the following, we present three general purpose preconditioners for solving the monolithic linear system~\eqref{eq:linearized-coupled-problem}.   {\color{rev1} The presented methods are preconditioners presented previously in the literature with the common feature that they are able to be implemented for the generic coupled system \eqref{eq:linearized-coupled-problem}.} For the sake of simplicity, system \eqref{eq:linearized-coupled-problem} is written here as a generic system of equations with block structure:
\begin{equation}
\label{eq:the-blocked-system}
\underbrace{\left[
\begin{array}{ccc}
\Amat_{11} &  \cdots  & \Amat_{1N} \\
\vdots     &  \ddots  & \vdots     \\
\Amat_{N1} &  \cdots  & \Amat_{NN}
\end{array}
\right]
}_{\displaystyle\Amat}
\underbrace{
\left[
\begin{array}{c}
\xmat_1\\ \vdots \\ \xmat_{N}
\end{array}
\right]
}_{\displaystyle\xmat}
=
\underbrace{
\left[
\begin{array}{c}
\bmat_1\\ \vdots \\ \bmat_{N}
\end{array}
\right]
}_{\displaystyle\bmat}.
\end{equation}
By setting $\Amat_{ij}:=\partial\fmat_i/\partial\umat_j$, $\xmat_i:=\Delta\umat_i$, and $\bmat_i:=-\fmat_i$, equation \eqref{eq:linearized-coupled-problem} is recovered. 

\subsection{Block Gauss-Seidel (BGS) preconditioner}
\label{sec:BGS}

The first preconditioner is based on a standard BGS iteration for solving the system~\eqref{eq:the-blocked-system}:
\begin{equation}
\label{eq:stationary:method:BGS}
\xmat^{(k)} = \xmat^{(k-1)} + \Mmat^{-1}_\mathrm{BGS}(\bmat -\Amat\xmat^{(k-1)}), \quad k=1,\ldots,k_\mathrm{max},
\end{equation}
where $\xmat^{(k)}$ is the $k$-th iterate in the solution process,  $\xmat^{(0)}$ is the initial guess  and $k_\mathrm{max}$ is the selected number of iterations. The (forward) BGS iteration is recovered by setting $\Mmat_\mathrm{BGS}$ as: 
\begin{equation*}
\Mmat_\mathrm{BGS} :=
\left[
\begin{array}{ccc}
\Amat_{11} & \cdots   & \zeromat \\
 \vdots    &  \ddots  &    \vdots  \\
\Amat_{N1} &  \cdots  & \Amat_{NN}
\end{array}
\right].
\end{equation*}
Matrix $\Mmat_\mathrm{BGS}$ is an approximation of matrix $\Amat$ obtained by discarding  the blocks placed at the upper triangle. The same idea applies by discarding the lower triangle instead of the upper (referred as backward BGS) or alternating the lower and upper triangles at each iteration (referred to as symmetric BGS). In our case, the three versions are implemented so that the user can choose the best variant for the particular application. {\color{rev1} The quality of the BGS preconditioner strongly depends in the error committed when discarding the upper (or lower) off-diagonal blocks of matrix $\Amat$. Thus, the BGS method requires that the diagonal block are dominant, see e.g. \cite{AxelssonNeytcheva2013}.}

\begin{figure}[ht!]
\centering
\includegraphics[width=0.99\textwidth]{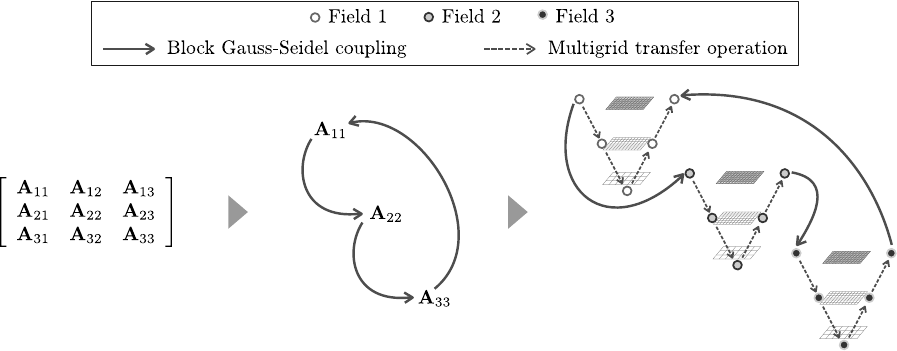}
\caption{Illustration of the BGS(AMG) preconditioner for a coupled problem with 3 fields. The method is built using to main steps: First a Block {\color{rev1}Gauss}-Seidel (BGS) iteration is considered for handling the coupling. Second, the resulting uncoupled systems are attacked with standard algebraic multigrid (AMG) methods.}
\label{fig:BGS-AMG}
\end{figure}

The forward BGS iteration \eqref{eq:stationary:method:BGS} is written in more detail as
\begin{equation}
\label{eq:bgs:iteration}
\xmat^{(k)}_i = \Amat^{-1}_{ii} \left[\bmat_i - \sum_{j=1}^{i-1} \Amat_{ij}\xmat^{(k)}_j + \sum_{j=i+1}^N \Amat_{ij}\xmat^{(k-1)}_j\right] \quad \text{for }i=1,\ldots,N.
\end{equation}
Thanks to the triangular structure of $\Mmat_\mathrm{BGS}$, computing the iteration \eqref{eq:bgs:iteration} only requires solving systems with the main diagonal blocks $\Amat_{ii}$, $i=1,\ldots,N$. That is, the BGS iteration transforms the coupled problem \eqref{eq:the-blocked-system} into a set of uncoupled single field-problems.  The advantage of this approach is that the resulting uncoupled systems can be solved with standard methods available for single field-problems, e.g. AMG methods tailored to each of the underlying fields. If AMG is considered for {\color{rev1} solving} the uncoupled problems, the resulting method is denoted as BGS(AMG), {\color{rev1} indicating} that the preconditioner is constructed using two main ingredients: an outer BGS iteration  plus inner uncoupled AMG methods (see Figure \ref{fig:BGS-AMG}).   
In our implementation, we allow the user to chose from several different methods for solving the single-field problems. The available methods include  aggregation-based AMG (provided by the MueLu package \citep{ProkopenkoHuWiesnerSiefertTuminaro2014}), relaxation methods and incomplete factorizations (both provided by the Ifpack package \citep{SalaHeroux2005}), and direct solvers (provided by the Amesos package \citep{SalaStanleyHeroux2006}). 

The BGS(AMG) is a general solver in the sense that it can be implemented for an arbitrary number of fields, and therefore, it can be potentially applied to a wide range of coupled problems.  However,  the method requires that all the main diagonal matrices $\Amat_{ii}$ are non-singular, which is not always fulfilled (e.g. in saddle point problems). 
Another potential drawback for some problems is that coupling is only achieved at the finest level, as it can easily be seen from Figure~\ref{fig:BGS-AMG}.

\subsection{SIMPLE preconditioner for an arbitrary number of fields}
\label{sec:SIMPLE}

The next preconditioner is based on an idea similar to the semi-implicit method for  pressure-linked equations (SIMPLE) \citep{CarettoGosmanPatankarSpalding1973} which was originally proposed for solving the saddle point structure of the Navier-Stokes equations.
This method is required 
in situations where the BGS(AMG) preconditioner fails (e.g. when one of the matrices $\Amat_{ii}$ is singular or zero).
The SIMPLE method is designed for problems with a $2\times 2$ block structure like:
\begin{equation}
\label{eq:the-blocked-system-2x2}
\left[
\begin{array}{cc}
{\Amat}_{11}   & {\Amat}_{12} \\
{\Amat}_{21}   & {\Amat}_{22}
\end{array}
\right]
\left[
\begin{array}{c}
\xmat_1 \\ \xmat_2
\end{array}
\right]
=
\left[
\begin{array}{c}
\bmat_1\\ \bmat_{2}
\end{array}
\right].
\end{equation}
The nice property of the method is that the main diagonal matrix $\Amat_{22}$ can be singular or even zero, as it turns out in saddle point problems, see \citep{ElmanHowleShadidShuttleworthTuminaro2008}.

The SIMPLE method for solving system \eqref{eq:the-blocked-system-2x2} is the solver characterized by the iteration
\begin{multline}
\label{eq:stationary:method:SIMPLE}
\def\arraystretch{1.5}
\left[
\begin{array}{c}
{\xmat}^{(k)}_{1} \\
{\xmat}^{(k)} _{2}
\end{array}
\right]
= 
\left[
\begin{array}{c}
{\xmat}^{(k-1)}_{1} \\
{\xmat}^{(k-1)} _{2}
\end{array}
\right]
+\\
\Mmat^{-1}_\mathrm{SIMPLE}\left(
\left[
\begin{array}{c}
{\bmat}_{1} \\
{\bmat}_{2}
\end{array}
\right]
 -
\left[
\begin{array}{cc}
{\Amat}_{11}   & {\Amat}_{12} \\
{\Amat}_{21}   & {\Amat}_{22}
\end{array}
\right] 
\left[
\begin{array}{c}
{\xmat}^{(k-1)}_{1} \\
{\xmat}^{(k-1)} _{2}
\end{array}
\right]
\right), \quad k=1,\ldots,k_\mathrm{max},
\end{multline}
where the superscript $[\cdot]^{(k)}$ denotes the $k$-th iteration in the solution process and the matrix 
$\Mmat_\mathrm{SIMPLE}$ is the following block factorization of the system matrix $\Amat$:
\begin{equation}
\label{eq:stationary:Matrix:SIMPLE}
\Mmat_\mathrm{SIMPLE} := 
\begin{bmatrix}
\Amat_{11} & \zerovec\\  \Amat_{21} & \widetilde\Smat
\end{bmatrix}
\begin{bmatrix}
\Imat & (\widetilde{\Amat}_{11})^{-1}{\Amat}_{12}\\ \zerovec & \Imat
\end{bmatrix}.
\end{equation}
Here, matrix $\widetilde\Amat_{11}$ is an easy-to-invert approximation of $\Amat_{11}$ and $\widetilde\Smat=\Amat_{22} - \Amat_{21}(\widetilde\Amat_{11})^{-1}\Amat_{12}$ is the associated approximation of the Schur complement matrix. If $\widetilde\Amat_{11}=\Amat_{11}$, 
matrix $\Mmat_\mathrm{SIMPLE}$ coincides with the system matrix $\Amat$ and the method \eqref{eq:stationary:method:SIMPLE} converges in one iteration. Thus, the effectivity of the SIMPLE iteration depends on the quality of the approximation $\widetilde\Amat_{11}$. In the numerical examples, we construct $\widetilde{\Amat}_{11}$ following the so-called SIMPLEC \citep{DoormaalRaithby1984} approach which consists in taking  $\widetilde{\Amat}_{11}$ as the diagonal matrix whose entries are the absolute values of the row sums of ${\Amat}_{11}$.

\begin{figure}[ht!]
\centering
\includegraphics[width=0.99\textwidth]{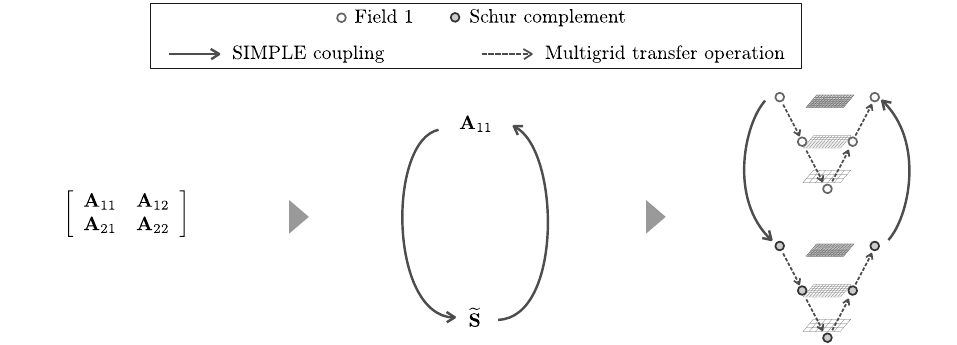}
\caption{Illustration of the SIMPLE(AMG) preconditioner for a coupled system with 2 fields. The method is obtained by applying an outer SIMPLE iteration for uncoupling the fields and AMG methods for solving the resulting uncoupled single-field systems.} 
\label{fig:simple-2x2}
\end{figure}

Iteration  \eqref{eq:stationary:method:SIMPLE} is written in more detail as
\begin{subequations}
\label{eq:SIMPLE}
\begin{eqnarray}
\Delta\xmat^{(k-1)}_1=& \Amat_{11}^{-1}\left[\bmat_1 - \Amat_{12}\xmat^{(k-1)}_2\right], &\text{(Predictor equation)}  \label{eq:SIMPLE:1}\\
\Delta\xmat^{(k-1)}_2 =&{\widetilde\Smat}^{-1}\left[\bmat_2-\Amat_{21}\Delta\xmat^{(k-1)}_1 -\Amat_{22}\xmat^{(k-1)}_2\right], &\text{(Schur equation)} \label{eq:SIMPLE:2}\\
\xmat^{(k)}_2 =& \xmat^{(k-1)}_2 + \Delta\xmat^{(k-1)}_2, &\text{(Update)}\\
\xmat^{(k)}_1 =& \xmat^{(k-1)}_1 - (\widetilde\Amat_{11})^{-1}\Amat_{12}\Delta\xmat^{(k-1)}_2. &\text{(Update)}
\end{eqnarray}
\end{subequations}
Steps \eqref{eq:SIMPLE:1} and \eqref{eq:SIMPLE:2} are known as the predictor and Schur equations and  require  solving systems with the matrices $\Amat_{11}$ and $\widetilde\Smat$ respectively.  Thus, the SIMPLE method can be seen as a procedure that transforms a coupled problem  \eqref{eq:the-blocked-system-2x2}  into a set of  two uncoupled systems that can be solved with standard techniques for single-field problems, e.g. AMG methods. The resulting preconditioner is denoted here as SIMPLE(AMG) highlighting that the methods consists in an outer SIMPLE iteration for uncoupling the fields and inner AMG solvers for the resulting uncoupled systems, see Figure \ref{fig:simple-2x2}. 

{\color{rev1} The SIMPLE method can be extended to cope} with generic coupled problems like \eqref{eq:the-blocked-system} with more than two physical fields (i.e. $N>2$). To this end, we follow a similar strategy to the one proposed by \citet{BadiaMartinPlanas2014} in the context of magneto-hydro dynamics. The methodology consists in merging the blocks of the system \eqref{eq:the-blocked-system} in order to recover a $2\times 2$ block structure. For instance, a  coupled problem with a $3\times 3$ block structure can be rearranged into a $2\times 2$ structure by merging the first and second fields:
\begin{equation}
\label{eq:3x3-system}
\left[
\begin{array}{cc|c}
{\Amat}_{11}   & {\Amat}_{12} & {\Amat}_{13} \\
{\Amat}_{21}   & {\Amat}_{22} & {\Amat}_{23} \\\hline
{\Amat}_{31}   & {\Amat}_{32} & {\Amat}_{33}
\end{array}
\right]
\left[
\begin{array}{c}
\xmat_1 \\ \xmat_2\\\hline \xmat_3
\end{array}
\right]
=
\left[
\begin{array}{c}
\bmat_1\\ \bmat_{2} \\\hline \bmat_3
\end{array}
\right]
\quad
\longrightarrow
\quad
\left[
\begin{array}{cc}
\overline{\Amat}_{11}   & \overline{\Amat}_{12} \\
\overline{\Amat}_{21}   & \overline{\Amat}_{22}
\end{array}
\right]
\left[
\begin{array}{c}
\overline\xmat_1 \\ \overline\xmat_2
\end{array}
\right]
=
\left[
\begin{array}{c}
\overline\bmat_1\\ \overline\bmat_{2}
\end{array}
\right],
\end{equation}
where
\begin{equation*}
\overline\Amat_{11} :=
\left[
\begin{array}{cc}
{\Amat}_{11}   & {\Amat}_{12} \\
{\Amat}_{21}   & {\Amat}_{22}
\end{array}
\right],\quad
\overline\Amat_{12}: =
\left[
\begin{array}{c}
{\Amat}_{13} \\
{\Amat}_{23}
\end{array}
\right],\quad
\overline\Amat_{21} :=
\left[
\begin{array}{cc}
{\Amat}_{31}   & {\Amat}_{32} \\
\end{array}
\right],\quad
\overline\Amat_{22} :=
\Amat_{33},
\end{equation*}
and 
\begin{equation*}
\overline\xmat_{1} :=
\left[
\begin{array}{c}
{\xmat}_{1} \\
{\xmat}_{2}
\end{array}
\right],\quad
\overline\xmat_{2} :=
\xmat_3,\quad
\overline\bmat_{1} :=
\left[
\begin{array}{c}
{\bmat}_{1} \\
{\bmat}_{2}
\end{array}
\right],\quad
\overline\bmat_{2} :=
\bmat_{3}.
\end{equation*}

%
%

\begin{figure}[ht!]
\centering
\includegraphics[width=0.99\textwidth]{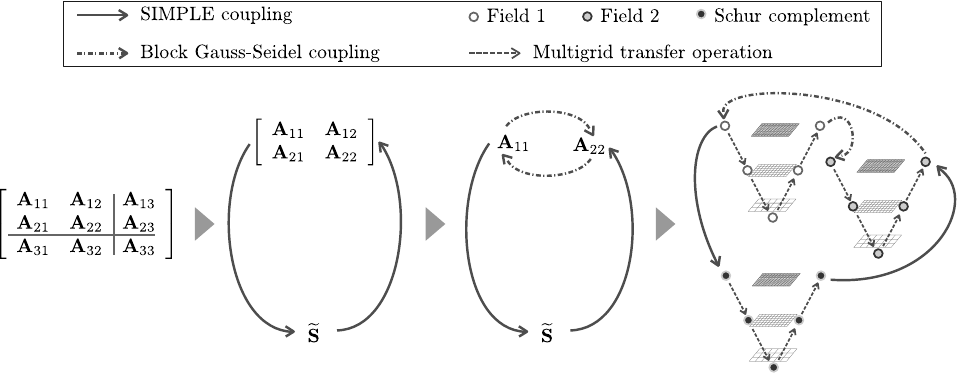}
\caption{Using the SIMPLE method in problems with more than 2 physical fields. The blocks of the original system are merged in order to recover the $2\times 2$ block structure required by the SIMPLE method. This results in a coupled problem to be solved within each SIMPLE iteration. A block iteration (e.g. BGS) is considered for uncoupling this problem, which results in a set of uncoupled single-field problems that can be efficiently solved with AMG.} 
\label{fig:simple-coupling-3x3}
\end{figure}

The merging procedure has to be performed carefully, because  after applying the SIMPLE iteration \eqref{eq:SIMPLE}, the resulting predictor and Schur equations \eqref{eq:SIMPLE:1} and \eqref{eq:SIMPLE:2} might involve more than one physical fields, and therefore standard AMG cannot be applied directly.
For example, the merging procedure for the $3\times 3$ problem given in equation \eqref{eq:3x3-system} leads to a system matrix for the predictor equation, which involves two physical fields:
$$\overline\Amat_{11} =
\left[
\begin{array}{cc}
{\Amat}_{11}   & {\Amat}_{12} \\
{\Amat}_{21}   & {\Amat}_{22}
\end{array}
\right].$$
Thus, in order to solve the system involving this matrix we cannot directly use standard AMG. Instead,  we consider a block iterative method,
either   BGS or SIMPLE, for uncoupling the fields. This can be seen as a recursive application of the different proposed approaches, which results in  uncoupled systems that can be solved with standard AMG methods, see Figure~\ref{fig:simple-coupling-3x3}.

\subsection{Monolithic AMG preconditioner}
\label{sec:AMG-BGS}

The third general purpose {\color{rev1} method} presented here is based on the so-called  monolithic AMG preconditioner proposed by  \citet{GeeKuettlerWall2011} for fluid-structure interaction (FSI) applications (see also \citet{MayrKloeppelWallGee2015}). Note that both the BGS(AMG) and SIMPLE(AMG) methods described above transform the coupled problem in a set of uncoupled systems which are {\color{rev1} solved} with independent AMG {\color{rev1} methods}.
The major drawback of this approach is that the coupling between the fields  is only resolved on the finest multigrid level through the outer block iteration. Even if efficient inner AMG solvers are used, the coupling can be under-resolved if not enough iterations of the outer BGS or SIMPLE are carried out.
This drawback is overcome by \citet{GeeKuettlerWall2011} in the context of FSI with a multigrid method which accounts for the coupling at all multigrid levels, while keeping the block structure of the fine level matrix. The monolithic AMG method is extended here for a generic coupled problem like~\eqref{eq:the-blocked-system}.

The main ingredients of the monolithic AMG method are standard AMG solvers for each one of the underlying physical fields.
Let us assume that we can build specific AMG methods for each one of the fields {\color{rev1} using information} of the  matrices $\Amat_{ii}$, $i=1,\ldots,N$. 
From the computed AMG hierarchies, we extract the usual multigrid objects:
1) the projection operators $\Pmat^\ell_i$,  2) the restriction operators $\Rmat^\ell_i$, and 3) the smoothers $\Svec^\ell_i$. The subscript $[\cdot]_i$ stands for each one of the fields $i=1,\ldots,N$, and the superscript $[\cdot]^\ell$ stands for each one of the levels $\ell=1,\ldots,\ell_\mathrm{max}$.
For the sake of simplicity, we assume that all the number of levels  $\ell_\mathrm{max}$ is  the same across all the fields.

The key idea of monolithic AMG is to reuse these single-field multigrid objects to build a multigrid method for the entire coupled problem. 
The transfer operators  $\Pmat^\ell$ and $\Rmat^\ell$  for the monolithic AMG method are defined as  the block diagonal matrices 
\begin{equation}
\label{eq:transfers-monolithic-AMG}
\Pmat^\ell := 
\left[
\begin{array}{ccc}
\Pmat^\ell_1 \\
 & \ddots & \\
 &  & \Pmat^\ell_N
\end{array}
\right],\quad
\Rmat^\ell := 
\left[
\begin{array}{ccc}
\Rmat^\ell_1 \\
 & \ddots & \\
 &  & \Rmat^\ell_N
\end{array}
\right],\quad
\ell=1,\ldots,\ell_\mathrm{max}-1,
\end{equation}
which can be computed from the transfer operators  $\Pmat^\ell_i$ and $\Rmat^\ell_i$ previously computed for the underlying single-field problems.

\begin{figure}[ht!]
\centering
\includegraphics[width=0.9\textwidth]{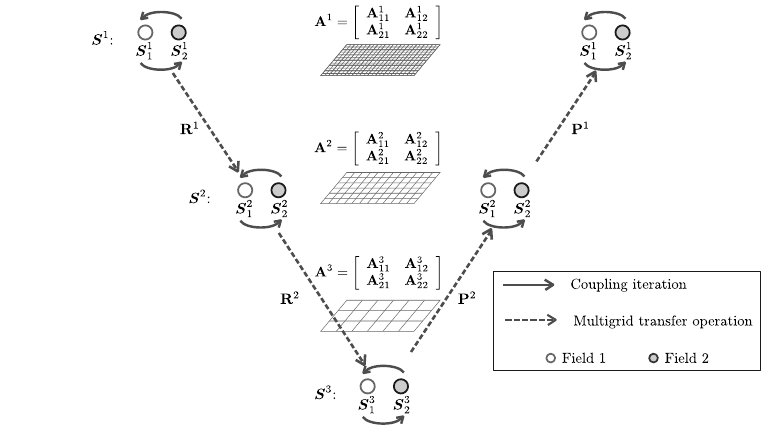}
\caption{Illustration of a monolithic AMG method for a coupled problem with 2 fields.} 
\label{fig:V-cycle-monolithic}
\end{figure}

The transfer operators $\Pmat^\ell$ and $\Rmat^\ell$  for the entire coupled problem have the nice property that lead to coarser system matrices $\Amat^\ell$ having the same block structure as the fine level matrix $\Amat^1$.
The coarse level matrices are obtained as
\begin{equation}
\underbrace{
\left[
\begin{array}{ccc}
\Amat^{\ell+1}_{11} &  \cdots  & \Amat^{\ell+1}_{1N} \\
\vdots     &  \ddots  & \vdots     \\
\Amat^{\ell+1}_{N1} &  \cdots  & \Amat^{\ell+1}_{NN}
\end{array}
\right]
}_{\displaystyle\Amat^{\ell+1}}
:=
\underbrace{
\left[
\begin{array}{ccc}
\Rmat^\ell_1 \\
 & \ddots & \\
 &  & \Rmat^\ell_N
\end{array}
\right]
}_{\displaystyle \Rmat^\ell}
\underbrace{
\left[
\begin{array}{ccc}
\Amat^\ell_{11} &  \cdots  & \Amat^\ell_{1N} \\
\vdots     &  \ddots  & \vdots     \\
\Amat^\ell_{N1} &  \cdots  & \Amat^\ell_{NN}
\end{array}
\right]
}_{\displaystyle\Amat^\ell}
\underbrace{
\left[
\begin{array}{ccc}
\Pmat^\ell_1 \\
 & \ddots & \\
 &  & \Pmat^\ell_N
\end{array}
\right],
}_{\displaystyle\Pmat^\ell}
\end{equation}
$\ell=1,\ldots,\ell_\mathrm{max}-1$, starting from the fine level matrix $\Amat^1$. 
Thanks to the diagonal structure of operators $\Pmat^\ell$ and $\Rmat^\ell$, the blocks of the coarse level matrices $\Amat^{\ell+1}$   are computed form the blocks of the fine level matrix $\Amat^\ell$, namely $\Amat^{\ell+1}_{ij}=\Rmat^\ell_i\Amat^\ell_{ij}\Pmat^\ell_j$, without mixing up the fields. That is, the physical interpretation of the coarse level blocks $\Amat^{\ell+1}_{ij}$ is the same as the fine level blocks $\Amat^\ell_{ij}$. {\color{rev1} This implies that all coarse level systems with matrices $\Amat^\ell$  have the same block structure as the fine level problem \eqref{eq:the-blocked-system}, see Figure~\ref{fig:V-cycle-monolithic}}.

Finally, the smoothers $\Smat^\ell$ used within the monolithic AMG method have to be defined. Since the coarse level matrices retain the block structure of the fine level problem, a block iterative method (i.e. BGS or SIMPLE) is used. More precisely, the smoothers $\Svec^\ell$ are build as an outer iteration (i.e. BGS or SIMPLE) for uncoupling the blocks in the system matrices $\Amat^\ell$, and the resulting single-field uncoupled problems are solved with the smothers $\Svec^\ell_i$  previously extracted from the single-field multigrid hierarchies, see Figure~\ref{fig:V-cycle-monolithic}.

The monolithic AMG method is finally constructed as a standard V-cycle using the projectors $\Pmat^\ell$, restrictors  $\Rmat^\ell$, coarse level matrices $\Amat^\ell$, and smoothers $\Smat^\ell$, see Figure \ref{fig:V-cycle-monolithic}. This method is {\color{rev1} referred to as} AMG(BGS) or AMG(SIMPLE) denoting that the block iteration (BGS or SIMPLE) is performed within each level of the AMG hierarchy.

\section{Implementation remarks}
\label{sec:implementation}

The {\color{rev1} presented} preconditioners for coupled problems (i.e. BGS, SIMPLE, and monolithic AMG) are implemented in such a way that they can be combined in order to build other more complex preconditioners allowing to solve challenging applications. 
For instance, any of these three methods can be used within the SIMPLE iteration \eqref{eq:SIMPLE} to solve the predictor equation \eqref{eq:SIMPLE:1} or the Schur complement equation \eqref{eq:SIMPLE:1} if the problems to be solved involve more than one physical field. This can be seen as a recursive application of a block preconditioner within another block preconditioner requiring a flexible implementation. 
We benefit from this flexibility  later in Section \ref{sec:lung:theory} for building preconditioners for the lung model, where we consider both the BGS method and the monolithic AMG preconditioner  within a SIMPLE iteration.

\begin{figure}[ht!]
\centering
\includegraphics[width=0.99\textwidth]{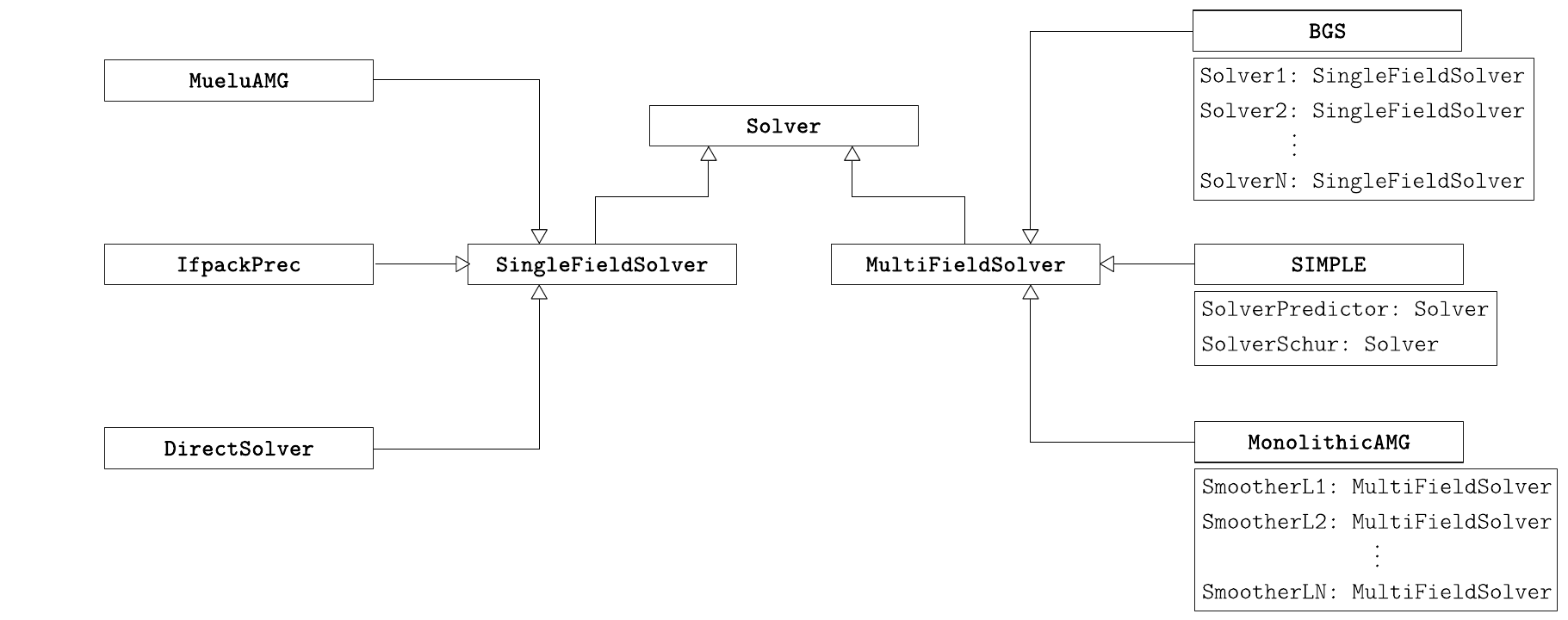}
\caption{ Simplified class diagram of the object oriented design of the proposed framework.}
\label{fig:class_diagram}
\end{figure}

The flexible implementation is mainly achieved considering an object-oriented software design, see Figure \ref{fig:class_diagram}. The framework is implemented using an interface class from which the BGS, SIMPLE, and monolithic AMG methods are derived. The solvers used within the SIMPLE iteration \eqref{eq:SIMPLE} to solve the predictor equation \eqref{eq:SIMPLE:1} or the Schur complement equation \eqref{eq:SIMPLE:1} are declared using the interface class, which allows calling BGS, SIMPLE, or monolithic AMG methods within the SIMPLE iteration. Similarly, the level smoothers within the monolithic AMG method are also declared using an interface class, allowing to have BGS or SIMPLE as smoothers. 
From the user point of view, the preconditioners are selected using a high level interface without the need of being familiar with the technical details of the underlying  implementation. 

\section{Applications}

In this section, we describe the different application examples that we have chosen to demonstrate the applicability and quality of the proposed approach.

\subsection{Thermo-structure interaction (TSI)}
\label{sec:TSI-theory}

The first problem type in which the proposed methods are tested is the interaction between an elastic solid and a thermal field. For the sake of a simple presentation, the solid is modeled here as a linear elastic material under small deformations, but the same approach holds for more sophisticated structural models. The solid is governed by the equations of elastodynamics
\begin{equation}
\label{eq:TSI:Solid}
\rho\ddot\dvec - \nabla\cdot \sigmavec(\dvec,u) = \bvec \quad\Omega\times I,
\end{equation}
describing the mechanical equilibrium in the domain $\Omega$ over the time interval $I$. The unknown displacements and accelerations are  $\dvec$ and $\ddot\dvec$ respectively, $\sigmavec$ stands for the unknown stresses, $\rho$ is the density, and $\bvec$ is an arbitrary body force.

Equation \eqref{eq:TSI:Solid} is augmented with the constitutive model
\begin{equation}
\label{eq:TSI:sterss}
\sigmavec(\dvec,u):=\boldsymbol{\mathcal{C}}:\varepsilonvec(\dvec) + m (u - u_0) \Ivec,
\end{equation}
where the term $\boldsymbol{\mathcal{C}}:\varepsilonvec(\dvec)$ is the standard elastic stress associated with the strains $\varepsilonvec(\dvec)$ via the fourth order elasticity tensor $\boldsymbol{\mathcal{C}}$, and the term $m (u - u_0) \Ivec$ is the thermal stress introducing the influence of the temperature field.
The value $u_0$ is the initial temperature, $\Ivec$ is the second order {\color{rev1}identity} tensor, and $m$ is defined as $m:=-(1\mu + 3\lambda)\alpha$, where $\lambda$ and $\mu$ are the Lam\'e constants, and $\alpha$ is the so-called \emph{coefficient of thermal expansion}. 

The TSI problem is obtained by coupling the mechanical equations \eqref{eq:TSI:Solid} and  \eqref{eq:TSI:sterss} with the time-dependent heat equation:
\begin{equation}
\label{eq:TSI:Thermo}
\rho C \dot u  - \nabla\cdot\qvec(u) - r(\dot\dvec)u = s \quad \Omega\times I,
\end{equation}
where $C$ is the heat capacity,  $\qvec(u)=-\kmat\nabla u$ is the heat flux, $\kmat$ is the second order tensor describing the thermal conductivity, and $s$ is an arbitrary source term. The influence of the mechanical displacements in the heat equation is modeled by the deformation-dependent reaction coefficient
\begin{equation*}
r(\dot\dvec) := m \left(\Ivec:\varepsilonvec(\dot\dvec)\right).
\end{equation*}
For the sake of brevity, it is assumed that the structural equation \eqref{eq:TSI:Solid} and the heat equation \eqref{eq:TSI:Thermo}  are closed with the usual boundary and initial conditions.

Equations \eqref{eq:TSI:Solid} and \eqref{eq:TSI:Thermo} lead to a fully coupled problem: the structural equation depends on the temperature and vice versa. Several solution strategies are proposed in the literature for resolving this coupling. They include staggered \citep{ArmeroSimo1992,ErbtsDuester2012,FarhatParkDubois-Pelerin1991,Miehe1995243} and monolithic methods \citep{DanowskiGravemeierYoshiharaWall2013}. In strongly coupled settings (i.e. under strong values of the thermal expansion coefficient $\alpha$), the work by \citet{DanowskiGravemeierYoshiharaWall2013} shows that monolithic methods {\color{rev1} are often superior to} partitioned schemes.

The monolithic scheme introduced by \citet{DanowskiGravemeierYoshiharaWall2013} leads to a system of equations with block structure requiring special preconditioners as the ones previously presented. The derivation of this linear system of equations is briefly presented here, see \citep{DanowskiGravemeierYoshiharaWall2013} for further details.

The discretization in space and time of 
equations \eqref{eq:TSI:Solid} and \eqref{eq:TSI:Thermo} lead to a set of non-linear equations that have to be solved at each time step in order to obtain the displacements  and temperatures. These equations are denoted as
\begin{subequations}
\label{eq:TSI:Nonlineqs}
\begin{eqnarray}
\fmat_\mathrm{S}^n(\dmat^n,\umat^n)=\zeromat,\\
\fmat_\mathrm{T}^n(\dmat^n,\umat^n)=\zeromat,
\end{eqnarray}
\end{subequations}
where $\fmat_\mathrm{S}^n$ and $\fmat_\mathrm{T}^n$ are the discrete structural and thermal equations, and $\dmat^n$ and $\umat^n$ are the discrete displacements and temperatures at a generic time point $t_n$.
Equations \eqref{eq:TSI:Nonlineqs} are a discrete, fully coupled, non-linear problem which can be written more compactly as
\begin{equation*}
\fmat^\mathrm{TSI}(\xmat^\mathrm{TSI})
=\zeromat,
\quad\text{with}\quad
\fmat^\mathrm{TSI}=
\begin{bmatrix}
\fmat_\mathrm{S}\\\fmat_\mathrm{T} 
\end{bmatrix}
\text{ and }
\xmat^\mathrm{TSI}=
\begin{bmatrix}
\dmat\\\umat 
\end{bmatrix},
\end{equation*}
where the  superscript $[\cdot]^n$ denoting the time point is dropped for simplicity. Following the monolithic approach by \citet{DanowskiGravemeierYoshiharaWall2013}, this non-linear problem is solved using a  Newton-type method
\begin{equation}
\label{eq:TSI:Newton}
\xmat^\mathrm{TSI}_{(i+1)}=\xmat^\mathrm{TSI}_{(i)} + \Delta\xmat^\mathrm{TSI}_{(i)}
\quad \text{with}\quad 
\left[\Partial{\fmat^\mathrm{TSI}}{\xmat^\mathrm{TSI}}\right]_{(i)}\Delta\xmat^\mathrm{TSI}_{(i)} = -\fmat^\mathrm{TSI}_{(i)},
\end{equation}
where the subscript $[\cdot]_{(i)}$ denotes the Newton iteration. The system of equations in \eqref{eq:TSI:Newton} involves the 
linearization of the TSI equations ($\fmat_\mathrm{S}$ and $\fmat_\mathrm{T}$) with respect the TSI variables ($\dmat$ and $\umat$). Thus, the system of linear equations to be solved at each Newton step has a 2$\times$2 block structure associated with the solid and thermal fields:
\begin{equation}
\label{eq:TSI:Monolithic}
\def\arraystretch{2.3}
\left[
\begin{array}{cc}
\Partial{\fmat_\mathrm{S}}{\dmat} & \Partial{\fmat_\mathrm{S}}{\umat}\\
\Partial{\fmat_\mathrm{T}}{\dmat} & \Partial{\fmat_\mathrm{T}}{\umat}
\end{array}
\right]_{(i)}
\left[
\begin{array}{c}
\Delta\dmat\\
\Delta\umat
\end{array}
\right]_{(i)}
=-
\left[
\begin{array}{c}
\fmat_\mathrm{S}\\\fmat_\mathrm{T} 
\end{array}
\right]_{(i)}.
\end{equation}
The main diagonal blocks of the system are the usual matrices arising in structural mechanics and in the heat equation, and therefore, they are available in many finite element packages. The off-diagonal blocks are associated with the coupling terms in equations \eqref{eq:TSI:Solid} and \eqref{eq:TSI:Thermo} and are computed as indicated in  \citep{DanowskiGravemeierYoshiharaWall2013}.

\begin{table}[ht!]
\centering
\begin{tabular}{lp{12cm}}
\toprule
Name & Description of the preconditioner\\
\midrule
BGS(AMG) &  Preconditioner described in Section \ref{sec:BGS}, using a backward BGS iteration.\\
SIMPLE(AMG) & Preconditioner described in Section \ref{sec:SIMPLE}, taking the structural equations as the predictor field and the thermal equations as Schur field. \\
AMG(BGS) & Preconditioner described in Section \ref{sec:AMG-BGS} with backward BGS smoothers.\\
\bottomrule
\end{tabular}
\caption{Description of the preconditioners for the monolithic thermo-structure interaction problem.}
\label{table:preconds:tsi}
\end{table}

The monolithic system \eqref{eq:TSI:Monolithic} {\color{rev1} requires often to be preconditioned in order to be solved with a Krylov method efficiently}. Since this system is a particular case of the general coupled problem \eqref{eq:the-blocked-system}, we can use the generic preconditioners previously presented, see Table~\ref{table:preconds:tsi}. The first method in Table \ref{table:preconds:tsi} is the BGS(AMG) preconditioner presented in Section \ref{sec:BGS} using a backwards BGS iteration. We choose a backward BGS since the effect of the structural displacement on the thermal equations is less important than the opposite effect, see~\citep{DanowskiGravemeierYoshiharaWall2013}. 
The second method is the SIMPLE(AMG) preconditioner introduced in Section \ref{sec:SIMPLE}. 
Finally, the third preconditioner is the monolithic AMG method presented in Section \ref{sec:AMG-BGS}. We consider a backwards BGS method as smoother within the AMG hierarchy and for that reason the method is denoted as AMG(BGS).

Among the three preconditioners given in Table~\ref{table:preconds:tsi}, the method BGS(AMG) corresponds to the one considered in the work by \citet{DanowskiGravemeierYoshiharaWall2013}. To the knowledge of the authors, the other two techniques are not previously used in the context of TSI.

\subsection{Fluid-structure interaction (FSI)}
\label{sec:Intro:FSI}

The second coupled problem considered in this work is FSI. We adopt the monolithic approach proposed by \citet{KuettlerWall2008} which furnishes challenging linear systems with block structure. For the sake of completeness, we briefly summarize the method here.

Let us consider the coupling between a solid and a fluid at the interface $\Gamma$. The displacements $\dvec^\mathrm{S}$ of the solid are governed by the equations of structural mechanics in the solid domain $\Omega^\mathrm{S}$. On the other hand,  the fluid velocities $\uvec^\mathrm{F}$ and pressure $p$ are governed by the arbitrary Lagrangian-Eulerian (ALE) version of the Naiver-Stokes equations \citep{FoersterWallRamm2006} in the fluid domain $\Omega^\mathrm{F}$. The ALE approach requires introducing a third (non-physical) grid movement field $\dvec^\mathrm{G}$  in $\Omega^\mathrm{F}$ in order to accommodate the fluid mesh to deformations of the FSI interface. The mesh movement $\dvec^\mathrm{G}$ is an arbitrary but unique extension of the solid displacement into the fluid domain $\Omega^\mathrm{F}$.

After discretizing in space and time, the FSI problem reduces to find the discrete structural displacements $\dmat^\mathrm{S}_n$ 
in $\Omega^\mathrm{S}$, the discrete grid displacements $\dmat^\mathrm{G}_n$ 
in $\Omega^\mathrm{F}$, and the discrete fluid velocities $\umat^\mathrm{F}_n$  and discrete fluid pressures $\pmat^\mathrm{F}_n$  
in $\Omega^\mathrm{F}$, for each time point $t_n$.

These unknowns are determined by solving the non-linear equations resulting from the discretization of the structural, fluid and grid motion fields. These equations are denoted here as:
\begin{subequations}
\label{eq:FSI:Interior}
\begin{eqnarray}
\fmat^\mathrm{S}_{\mathrm{I},n}(\dmat^\mathrm{S}_n) = \zeromat,\\
\fmat^\mathrm{G}_{\mathrm{I},n}(\dmat^\mathrm{G}_n) = \zeromat,\\
\fmat^\mathrm{F}_{\mathrm{I},n}(\umat^\mathrm{F}_n,\pmat^\mathrm{F}_n,\dmat^\mathrm{G}_n) = \zeromat,
\end{eqnarray}
\end{subequations}
where the superscripts $[\cdot]^\mathrm{S}$, $[\cdot]^\mathrm{G}$ and $[\cdot]^\mathrm{F}$ stand for the structural field, the grid motion and the fluid  respectively. The subscript $[\cdot]_\mathrm{I}$ denotes that the equations are associated only with mesh nodes in the interior of the respective domains (i.e. not at the FSI interface).

Equations \eqref{eq:FSI:Interior} have to be augmented with boundary and initial conditions together with the FSI conditions at the interface $\Gamma$ in order to have a well posed problem. Boundary and initial conditions are standard, and therefore, they are not discussed here.
The coupling conditions at the FSI interface $\Gamma$ are:
\begin{subequations}
\label{eq:FSI:FSIcondCont}
\begin{eqnarray}
\uvec^\mathrm{F} = \dot{\dvec}^\mathrm{S}\quad\text{on }\Gamma\times I,\label{eq:FSI:FSIcondCont:1}\\
\dvec^\mathrm{G} = \dvec^\mathrm{S}\quad\text{on }\Gamma\times I,\label{eq:FSI:FSIcondCont:2}\\
\sigmavec^\mathrm{S}\cdot\nvec =  \sigmavec^\mathrm{F}\cdot\nvec \quad\text{on }\Gamma\times I.\label{eq:FSI:FSIcondCont:3}
\end{eqnarray}
\end{subequations}
The first condition is the kinematic continuity between fluid and structural velocities, the second is the continuity between structural displacements and the grid motion field, and the third is the equilibrium of forces at the FSI interface. The tensors $\sigmavec^\mathrm{S}$ and $\sigmavec^\mathrm{F}$ stand for the Cauchy stress of the solid and the fluid respectively, and $\nvec$ is a unique normal defined on the FSI interface. 

{\color{rev1} Defining} the discrete FSI problem requires introducing the discrete version of the FSI coupling condition \eqref{eq:FSI:FSIcondCont:3} which is denoted here as
\begin{equation}
\label{eq:FSI:CouplDisc:forces}
\fmat^\mathrm{S}_{\mathrm{\Gamma},n}(\dmat^\mathrm{S}_n) + \fmat^\mathrm{F}_{\mathrm{\Gamma},n}(\umat^\mathrm{F}_n,\pmat^\mathrm{F}_n,\dmat^\mathrm{G}_n) = \zeromat,
\end{equation}
where the subscript $[\cdot]_\Gamma$ denotes quantities associated with the FSI interface $\Gamma$. The vector $\fmat^\mathrm{S}_{\mathrm{\Gamma},n}$ represents the structural interface forces obtained by weighting (the discrete version of) $\sigmavec^\mathrm{S}\cdot\nvec$ by the shape functions associated with the nodes on $\Gamma$. The same applies for the fluid interface forces $\fmat^\mathrm{F}_{\mathrm{\Gamma},n}$.

The interior field equations \eqref{eq:FSI:Interior} together with the interface condition \eqref{eq:FSI:CouplDisc:forces} form the fully coupled FSI problem to be solved at each time step:
\begin{equation}
\label{eq:FSI:nonlinproblem-4x4}
\def\arraystretch{1.4}
\fmat^\mathrm{FSI}(\xmat^\mathrm{FSI}) = \zeromat
\quad\text{ with }
\fmat^\mathrm{FSI} :=
\left[
\begin{array}{c}
\fmat^\mathrm{S}_\mathrm{I}\\
\fmat^\mathrm{G}_\mathrm{I}\\
\fmat^\mathrm{F}_\mathrm{I}\\
\fmat^\mathrm{S}_\mathrm{\Gamma} + \fmat^\mathrm{F}_\mathrm{\Gamma}
\end{array}
\right]
\text{ and }
\xmat^\mathrm{FSI} :=
\left[
\begin{array}{c}
\dmat^\mathrm{S}_\mathrm{I}\\ \dmat^\mathrm{G}_\mathrm{I} \\ \umat^\mathrm{F}_\mathrm{I}\\ \umat^\mathrm{F}_\mathrm{\Gamma}
\end{array}
\right],
\end{equation}
where the subscript $[\cdot]_n$ for the time step has been dropped for simplicity. Note that the FSI variable $\xmat^\mathrm{FSI}$ does not include the displacements $\dmat^\mathrm{S}_{\mathrm{\Gamma},n}$ nor $\dmat^\mathrm{G}_{\mathrm{\Gamma},n}$. This is because these displacements can be computed once $\umat^\mathrm{F}_{\mathrm{\Gamma},n}$ is available using the coupling conditions \eqref{eq:FSI:FSIcondCont:1} and \eqref{eq:FSI:FSIcondCont:2}. For the ease of notation, the fluid pressure $\pmat^\mathrm{F}$ is also dropped from the formula and it is assumed that the fluid quantities $\umat^\mathrm{F}_\mathrm{I}$ and $\umat^\mathrm{F}_\mathrm{\Gamma}$ include  the pressure degrees of freedom.

By merging the interior and interface fluid degrees of freedom, namely
\begin{equation*}
\def\arraystretch{1.4}
\umat^\mathrm{F}:=
\left[
\begin{array}{c}
\umat^\mathrm{F}_\mathrm{I}\\
\umat^\mathrm{F}_\mathrm{\Gamma}
\end{array}
\right]
\quad\text{and}\quad
\fmat^\mathrm{F}:=
\left[
\begin{array}{c}
\fmat^\mathrm{F}_\mathrm{I}\\
\fmat^\mathrm{S}_\mathrm{\Gamma} + \fmat^\mathrm{F}_\mathrm{\Gamma}
\end{array}
\right],
\end{equation*}
the non-linear FSI problem \eqref{eq:FSI:nonlinproblem-4x4} is rewritten as
\begin{equation}
\label{eq:FSI:nonlinproblem}
\def\arraystretch{1.4}
\fmat^\mathrm{FSI}(\xmat^\mathrm{FSI}) = \zeromat
\quad\text{ with }
\fmat^\mathrm{FSI} :=
\left[
\begin{array}{c}
\fmat^\mathrm{S}_\mathrm{I}\\
\fmat^\mathrm{G}_\mathrm{I}\\
\fmat^\mathrm{F}
\end{array}
\right]
\text{ and }
\xmat^\mathrm{FSI} :=
\left[
\begin{array}{c}
\dmat^\mathrm{S}_\mathrm{I}\\ \dmat^\mathrm{G}_\mathrm{I} \\ \umat^\mathrm{F}
\end{array}
\right].
\end{equation}

Following the approach by \citet{KuettlerWall2008}, the fully coupled FSI problem \eqref{eq:FSI:nonlinproblem} is solved using a Newton-like method:
\begin{equation}
\label{eq:FSI:Newton}
\xmat^\mathrm{FSI}_{(i+1)}=\xmat^\mathrm{FSI}_{(i)} + \Delta\xmat^\mathrm{FSI}_{(i)}
\quad \text{with}\quad 
\left[\Partial{\fmat^\mathrm{FSI}}{\xmat^\mathrm{FSI}}\right]_{(i)}\Delta\xmat^\mathrm{FSI}_{(i)} = -\fmat^\mathrm{FSI}_{(i)}.
\end{equation}
Note that the Jacobian matrix $\partial\fmat^\mathrm{FSI}/\partial\xmat^\mathrm{FSI}$ involves the linearization of $\fmat^\mathrm{S}_\mathrm{I}$,  $\fmat^\mathrm{G}_\mathrm{I}$, $\fmat^\mathrm{F}$  with respect to the variables $\dmat^\mathrm{S}_\mathrm{I}$, $\dmat^\mathrm{G}_\mathrm{I}$,  $\umat^\mathrm{U}$. This leads to a system of equations with the following block structure
\begin{equation}
\label{eq:FSI:BlockSystem}
\def\arraystretch{2.3}
\left[
\begin{array}{ccc}
\Partial{\fmat^\mathrm{S}_\mathrm{I}}{\dmat^\mathrm{S}_\mathrm{I}}   &  \zeromat & \Partial{\fmat^\mathrm{S}_\mathrm{I}}{\umat^\mathrm{U}}  \\
\zeromat & \Partial{\fmat^\mathrm{G}_\mathrm{I}}{\dmat^\mathrm{G}_\mathrm{I}} & \Partial{\fmat^\mathrm{G}_\mathrm{I}}{\umat^\mathrm{U}} \\
\Partial{\fmat^\mathrm{F}}{\dmat^\mathrm{S}_\mathrm{I}} & \Partial{\fmat^\mathrm{F}}{\dmat^\mathrm{G}_\mathrm{I}} & \Partial{\fmat^\mathrm{F}}{\umat^\mathrm{F}} \\
\end{array}
\right]
\left[
\begin{array}{c}
\Delta\dmat^\mathrm{S}_\mathrm{I}\\ \Delta\dmat^\mathrm{G}_\mathrm{I}\\ \Delta\umat^\mathrm{U}
\end{array}
\right]
=
- 
\left[
\begin{array}{c}
\fmat^\mathrm{S}_\mathrm{I}\\
\fmat^\mathrm{G}_\mathrm{I}\\
\fmat^\mathrm{F}
\end{array}
\right],
\end{equation}
corresponding to the solid, grid motion and fluid fields.

Solving the system of equations \eqref{eq:FSI:BlockSystem} is very challenging because it combines physical fields of very different kind. Particular preconditioners for this problem are proposed for example in \citep{CrosettoDeparisFouresteyQuarteroni2011,GeeKuettlerWall2011, Heil2004}. In the examples below, we compare the performance of the FSI-specific solvers by \citet{GeeKuettlerWall2011} with respect to the general purpose preconditioners and implementation proposed in Section \ref{sec:all-precs}. The considered methods are summarized in Table \ref{table:preconds:fsi}.

\begin{table}[ht!]
\centering
\begin{tabular}{lp{12cm}}
\toprule
Name & Description of the preconditioner\\
\midrule
BGS(AMG) generic &  General purpose preconditioner given in Section \ref{sec:BGS}.\\
AMG(BGS) generic & General purpose preconditioner given in Section \ref{sec:AMG-BGS} with BGS smoothers.\\
BGS(AMG) FSI & FSI-specific version of the preconditioner BGS(AMG) given in  \citep{GeeKuettlerWall2011}.\\
AMG(BGS) FSI & FSI-specific version of the preconditioner AMG(BGS) given in  \citep{GeeKuettlerWall2011}.\\
\bottomrule
\end{tabular}
\caption{Description of the preconditioners for fluid-structure interaction.}
\label{table:preconds:fsi}
\end{table}

\subsection{Lung model}
\label{sec:lung:theory}

The third and most challenging coupled problem considered in this paper is the lung model introduced by \citet{YoshiharaIsmailWall2013,YoshiharaRothWall2015}.
 This model is an extension of a standard FSI problem in order to account for the inflation of the lung tissue in the respiration process.  The key rationale of the model is the following.  
 
\begin{figure}[ht!]
\centering
\includegraphics[width=0.6\textwidth]{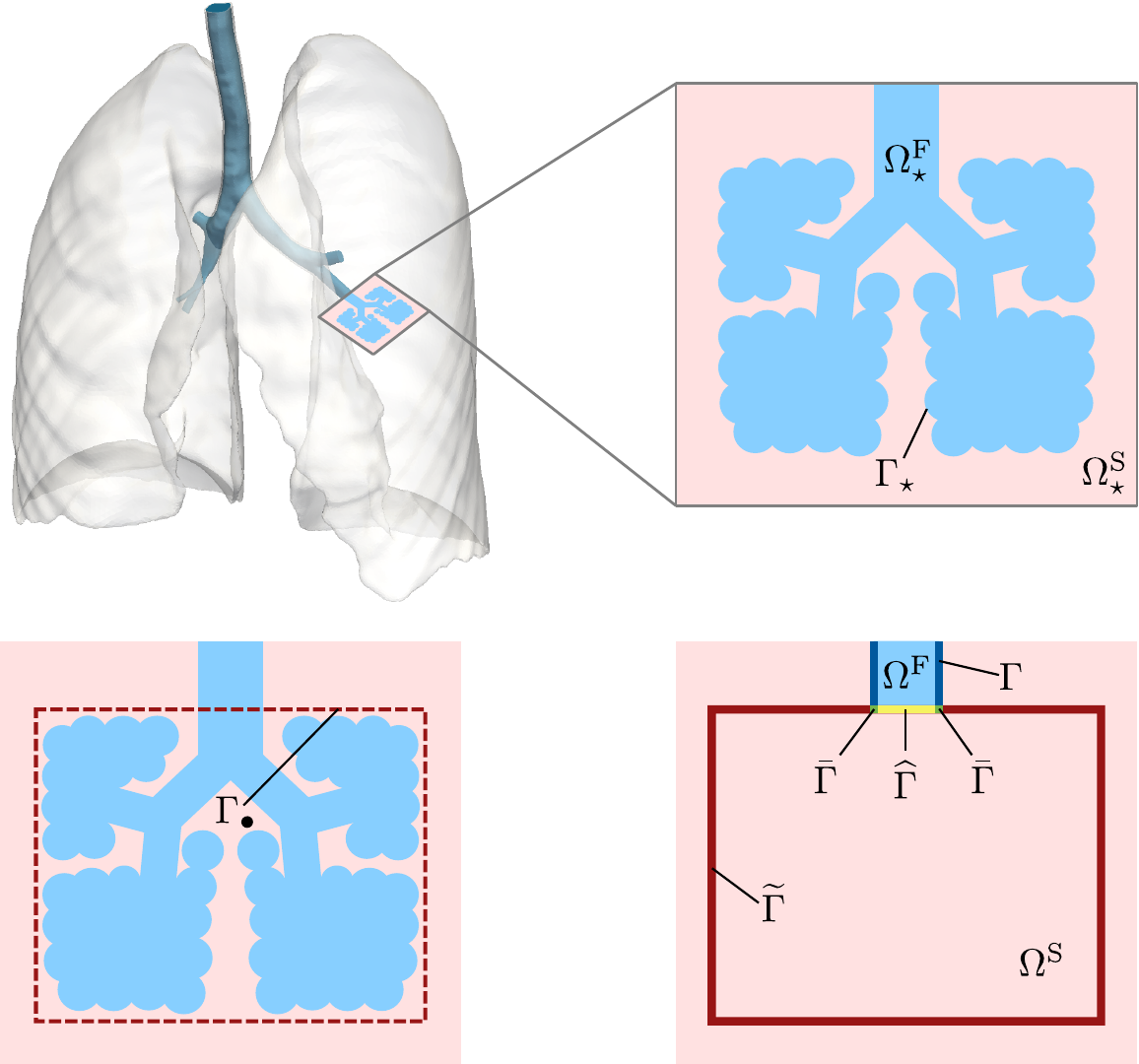}
\caption{Illustration of the lung model introduced by \citet{YoshiharaIsmailWall2013,YoshiharaRothWall2015}.}
\label{fig:Lung:Model}
\end{figure}

The human lung can be ideally modeled as an FSI problem  coupling the air transported by the respiratory tree and the lung tissue. Completely resolving this problem is usually not feasible due to the complex geometries involved: only few branches of the airway tree are imageable with CT-scans, and therefore, the domains $\Omega^\mathrm{S}_\star$, $\Omega^\mathrm{F}_\star$ and $\Gamma_\star$ of the ideal FSI problem are not available, see Figure \ref{fig:Lung:Model} (top). In the computational model, the airways end in artificial outlets, but in reality, the branching of the conducting passages continues until the alveolar region is reached. The behavior of this region completely dominates the overall behavior of the system.
Following the approach by \citet{YoshiharaIsmailWall2013,YoshiharaRothWall2015}, the unresolved structures are surrounded by an 
 artificial boundary $\Gamma_\bullet$, see Figure \ref{fig:Lung:Model} (bottom, left), and a suitable model  is introduced inside the delimited region accounting for the unresolved features.  The region enclosed by $\Gamma_\bullet$ is modeled as an elastic solid with the same properties as the surrounding lung tissue,  see Figure \ref{fig:Lung:Model} (bottom, right).
The discarded fluid domain is taken into account by assuming that the air transported in the feeding vessel ending at $\widehat\Gamma$  causes an inflation of the region delimited by  $\Gamma_\bullet$.
That is, the volume change of the region enclosed by  $\Gamma_\bullet$ is equal to the volume of air
flowing into it through $\widehat\Gamma$. This constraint is expressed formally as
\begin{equation}
\label{eq:FSI:constraintCont}
\dot V_{\Gamma_\bullet} = \int_{\widehat\Gamma} (\uvec^\mathrm{F}-\dot\dvec^\mathrm{G})\cdot\nvec\text{\ d}A,
\end{equation}
where $\dot V_{\Gamma_\bullet}$ denotes the rate of change of the volume surrounded by $\Gamma_\bullet$ and the right hand side is the flow rate through the outlet $\widehat\Gamma$. 

In order to have a well defined problem, extra conditions have to be introduced on the artificial ends of the airways:
\begin{subequations}
\label{eq:FSI:constraintExtraCond}
\begin{eqnarray}
\uvec^\mathrm{F} =& \dot\dvec^\mathrm{S}&\quad\text{on }\bar\Gamma,\\
\dvec^\mathrm{G} =& \dvec^\mathrm{S}&\quad\text{on }\bar\Gamma\cup\widehat\Gamma,\\
\sigmavec^\mathrm{F}\cdot\nvec =& \zerovec&\quad\text{on }\widehat\Gamma,\\
\sigmavec^\mathrm{S}\cdot\nvec =& \zerovec&\quad\text{on }\widehat\Gamma.
\end{eqnarray}
\end{subequations}
The first condition is the standard continuity of velocities {\color{rev1} inherited} from the FSI interface. The second condition is introduced in order to avoid gaps between the fluid and solid domains on the artificial boundary $\bar\Gamma\cup\widehat\Gamma$, and the two last are  standard ``do {\color{rev1} nothing}'' conditions.

In conclusion, the lung model results in an FSI problem which is augmented with the constraint \eqref{eq:FSI:constraintCont} and the extra conditions \eqref{eq:FSI:constraintExtraCond}. This leads (after the usual space and time discretization) to a non-linear system of algebraic equations to be solved at each time step. {\color{rev1} Defining} this problem requires introducing the discrete version of constraint \eqref{eq:FSI:constraintCont}, which is denoted here using the scalar equation
\begin{equation*}
c(\dmat^\mathrm{S},\umat^\mathrm{F})=0.
\end{equation*}
The number of such  constraints in the model is equal to the number of artificial outlets $\widehat\Gamma$. For convenience, the constraints corresponding to all the artificial outlets are grouped in the vector equation
\begin{equation}
\label{eq:FSI:discrete:cons}
\cmat(\dmat^\mathrm{S},\umat^\mathrm{F})=0.
\end{equation}

With {\color{rev1} this notation}, the non-linear algebraic equations to be solved at each time step are
\begin{equation*}
\fmat^\mathrm{C;FSI}(\xmat^\mathrm{C;FSI}) = \zeromat,
\end{equation*}
where
\begin{equation}
\label{eq:FSI:nonlinproblem4}
\def\arraystretch{2}
\fmat^\mathrm{C;FSI} =
\left[
\begin{array}{c}
\fmat^\mathrm{S}_\mathrm{I\cup\widehat\Gamma} \\
\fmat^\mathrm{G}_\mathrm{I}   \\
\fmat^\mathrm{F}_\mathrm{I\cup\widehat\Gamma}\\
\fmat^\mathrm{S}_\mathrm{\Gamma\cup\bar\Gamma} + \fmat^\mathrm{F}_\mathrm{\Gamma\cup\bar\Gamma}\\
\cmat
\end{array}
\right]+ 
\left[
\begin{array}{c}
\lambdavec^\mathrm{T}\Partial{\cmat}{\dmat^\mathrm{S}_\mathrm{I\cup\widehat\Gamma}}\\
 \lambdavec^\mathrm{T}\Partial{\cmat}{\dmat^\mathrm{G}_\mathrm{I}}  \\
\lambdavec^\mathrm{T} \Partial{\cmat}{\umat^\mathrm{F}_\mathrm{I\cup\widehat\Gamma}} \\
 \lambdavec^\mathrm{T}\Partial{\cmat}{\umat^\mathrm{F}_\mathrm{\Gamma\cup\bar\Gamma}} \\
\zerovec
\end{array}
\right]
\quad\text{ and }\quad
\xmat^\mathrm{C;FSI} =
\left[
\begin{array}{c}
\dmat^\mathrm{S}_\mathrm{I\cup\widehat\Gamma}\\
\dmat^\mathrm{G}_\mathrm{I} \\ \umat^\mathrm{F}_\mathrm{I\cup\widehat\Gamma}\\ \umat^\mathrm{F}_\mathrm{\Gamma\cup\bar\Gamma}\\ \lambdavec
\end{array}
\right].
\end{equation}
Note that the first four equations are essentially the standard FSI equations \eqref{eq:FSI:nonlinproblem-4x4} extended with the interface conditions \eqref{eq:FSI:constraintExtraCond} and  the loads introduced by the Lagrange  multipliers $\lambdavec$ associated  with the constraint \eqref{eq:FSI:discrete:cons}. The last  equation in \eqref{eq:FSI:nonlinproblem4} corresponds to the algebraic constraint  $\cmat$. After linearization of problem \eqref{eq:FSI:nonlinproblem4} within a Newton-like procedure, 
the system to be solved at each Newton step has the form:
\begin{equation}
\label{eq:FSI:BlockSystem:lung}
\def\arraystretch{1.4}
\left[
\begin{array}{ccc|c}
\Amat^\mathrm{SS}   &  \Amat^\mathrm{SG} & \Amat^\mathrm{SF} & \Amat^\mathrm{SC} \\
\Amat^\mathrm{GS}   &  \Amat^\mathrm{GG} & \Amat^\mathrm{GF} & \Amat^\mathrm{GC} \\
\Amat^\mathrm{FS}   &  \Amat^\mathrm{FG} & \Amat^\mathrm{FF} & \Amat^\mathrm{FC} \\
\hline
\Amat^\mathrm{CS}   &  \Amat^\mathrm{CG} & \Amat^\mathrm{CF} & \zerovec \\
\end{array}
\right]
\left[
\begin{array}{c}
\Delta\dmat^\mathrm{S}_\mathrm{I\cup\widehat\Gamma}\\ \Delta\dmat^\mathrm{G}_\mathrm{I}\\ \Delta\umat^\mathrm{F}\\ 
\hline
\Delta \lambdavec
\end{array}
\right]
=
- 
\left[
\begin{array}{c}
\bmat^\mathrm{S}\\
\bmat^\mathrm{G}\\
\bmat^\mathrm{F}\\
\hline \cvec
\end{array}
\right].
\end{equation}
This system of equations consists of 4 matrix equations associated with the 3 FSI fields (solid, grid motion and fluid) plus a block associated with the constraint. The FSI  blocks are essentially the same as  the ones already presented in the FSI system \eqref{eq:FSI:BlockSystem}. 
The precise {\color{rev1}definition} of all the objects in  \eqref{eq:FSI:BlockSystem:lung} is not provided here for brevity, refer to \citet{YoshiharaIsmailWall2013,YoshiharaRothWall2015} for details.

Solving the system of equations \eqref{eq:FSI:BlockSystem:lung} is very challenging. In the one hand, it shares the difficulties already commented for the pure FSI system \eqref{eq:FSI:BlockSystem}. Moreover, the system has a zero matrix on the main diagonal due to the saddle point structure introduced by the volumetric constraint \eqref{eq:FSI:discrete:cons}. This means that the BGS method cannot be used alone for solving the system, and therefore, the preconditioner BGS(AMG) introduced in Section \ref{sec:BGS} cannot be applied.

The preconditioners considered here for the lung model are based on the SIMPLE method (cf. Section~\ref{sec:SIMPLE}). Recall that in order to use  SIMPLE, the $4\times 4$ block structure of system \eqref{eq:FSI:BlockSystem:lung} has to be transformed into a $2\times2$ structure. This is achieved here by grouping the FSI blocks as indicated  in equation \eqref{eq:FSI:BlockSystem:lung}. The SIMPLE requires solving systems with the predictor matrix $\Amat_{11}$ and Schur complement matrix $\widetilde\Smat$.
For the splitting indicated  in equation \eqref{eq:FSI:BlockSystem:lung}, the matrix $\Amat_{11}$  is 
\begin{equation}
\label{eq:FSI-block}
\Amat_{11} = 
\left[
\begin{array}{ccc}
\Amat^\mathrm{SS}   &  \Amat^\mathrm{SG} & \Amat^\mathrm{SF}\\
\Amat^\mathrm{GS}   &  \Amat^\mathrm{GG} & \Amat^\mathrm{GF}\\
\Amat^\mathrm{FS}   &  \Amat^\mathrm{FG} & \Amat^\mathrm{FF}
\end{array}
\right],
\end{equation}
corresponding to the FSI operator. That is, at each iteration of the SIMPLE method, a linearized FSI problem has to be solved. To solve this problem we consider the two FSI preconditioners previously presented  in Table \ref{table:preconds:fsi}, namely BGS(AMG), and AMG(BGS).

On the other hand, the system involving the {\color{rev1} Schur} complement  $\widetilde\Smat$ is solved with a direct solver because the size of $\widetilde\Smat$  is frequently rather small in this application. This is because, the number of rows in $\widetilde\Smat$ corresponds to the number of ending outlets in the lung model which in practice is much smaller than the number of degrees of freedom in the other fields.

\begin{figure}[ht!]
\centering
\includegraphics[width=0.99\textwidth]{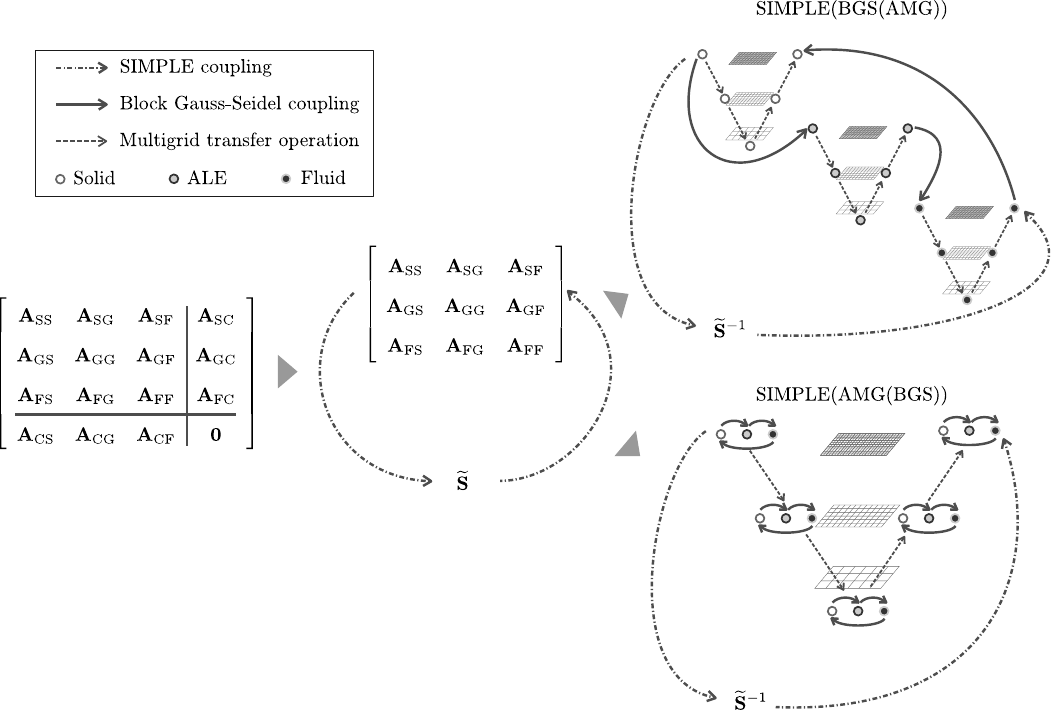}
\caption{Illustration of the preconditioners for the lung model. We consider two different methods, namely SIMPLE(BGS(AMG)) and SIMPLE(AMG(BGS). Both preconditioners are obtained by merging the fields corresponding to the FSI problem to recover a $2\times 2$ block structure in order to apply SIMPLE.  Within the SIMPLE iteration, a linearized FSI problem has to be solved. We consider two options, either a BGS(AMG) method enforcing the coupling at the finest multigrid level, or a AMG(BGS) method resolving the coupling at all levels.}
\label{fig:Lung:precs}
\end{figure}

In conclusion, we consider two different preconditioners for the lung model, see Table \ref{table:preconds:lung}. The first one is denoted as SIMPLE(BGS(AMG)) and consists in an outer SIMPLE iteration (Section~\ref{sec:SIMPLE}) and an inner BGS(AMG) preconditioner (cf. Section \ref{sec:BGS}) for solving the FSI block \eqref{eq:FSI-block}. The second preconditioner  is denoted as SIMPLE(AMG(BGS)) and is an outer SIMPLE iteration  and an inner monolithic AMG preconditioner (cf. Section \ref{sec:AMG-BGS}) for solving the FSI block. The only difference between these to methods it that SIMPLE(BGS(AMG)) resolves the FSI coupling only at the finest AMG level, while SIMPLE(AMG(BGS)) resolves the the FSI coupling at all levels, see Figure \ref{fig:Lung:precs}.

\begin{table}[ht!]
\centering
\begin{tabular}{lp{11.5cm}}
\toprule
Name & Description of the preconditioner\\
\midrule
SIMPLE(BGS(AMG)) & Outer SIMPLE iteration using an inner BGS(AMG) preconditioner  for the FSI subproblem.\\
SIMPLE(AMG(BGS)) & Outer SIMPLE iteration using an inner AMG(BGS) preconditioner for the FSI subproblem.\\
\bottomrule
\end{tabular}
\caption{Description of the preconditioners for the lung model.}
\label{table:preconds:lung}
\end{table}

%
%
%
%
%
%
%
%
%
%

\section{Numerical examples}

\subsection{TSI example: Second Danilovskaya problem}
\label{sec:TSI:example:danilovsakaya}

This first example is devoted to study the performance and scalability  of the {\color{rev1} proposed} preconditioners 
in a benchmark test for  TSI. The example is known as the \emph{second Danilovskaya problem}, it was originally proposed in \citep{Danilovskaya1952}, and it has often been used in literature for validating  fully coupled thermo-mechanical models (e.g. in \citep{DanowskiGravemeierYoshiharaWall2013,FarhatParkDubois-Pelerin1991,TanakaMatsumotoMoradi1995}).

\begin{figure}[ht!]
\centering
\includegraphics[width=0.9\textwidth]{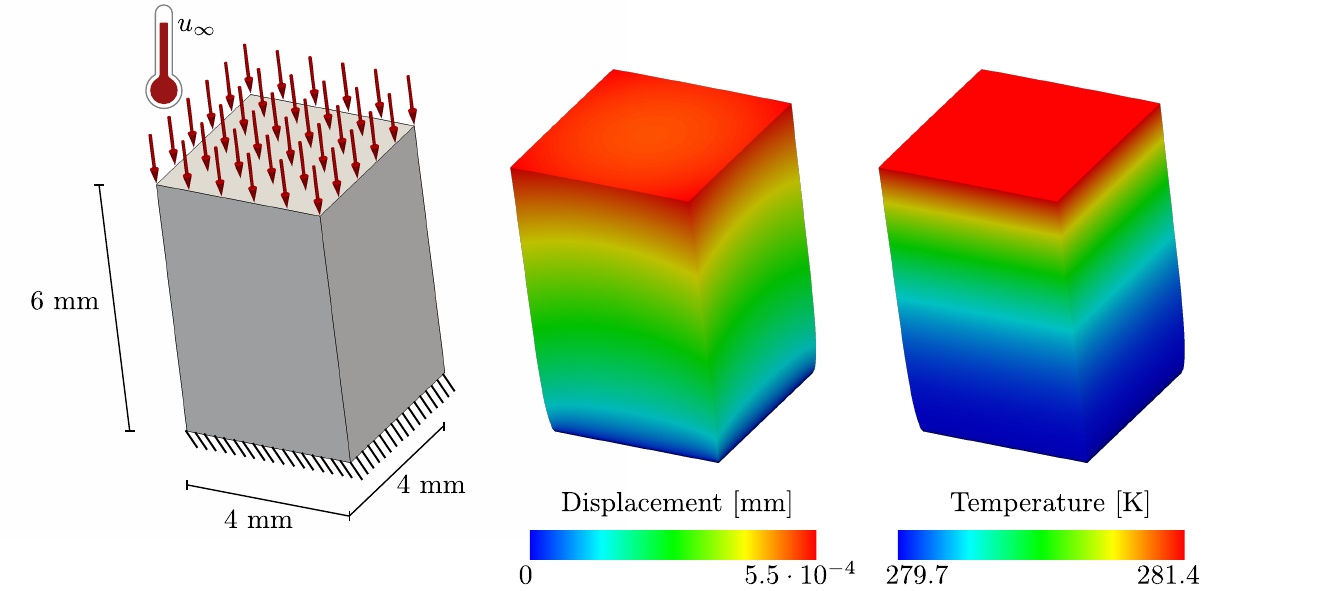}
\caption{Second Danilovskaya problem: Problem description (left), deformation (center) and temperature distribution (right) at time $t=4$ s.}
\label{fig:tsi-dani-statement}
\end{figure}

The example consists in a thermo-elastic  rectangular prism which is heated on the top face, see Figure~\ref{fig:tsi-dani-statement}. The thermo-mechanical properties of the material are Young's modulus $E=210\cdot 10^9$ Pa, Poisson's ratio $\nu=0.3$, density  $\rho=7.86\cdot 10^3$ kg/m$^3$, heat capacity $C=0.821$ J/(kg $\cdot$ K), heat conductivity $k=1.03$ {W/(m $\cdot$ K)} and coefficient of thermal expansion $\alpha = 1.1\cdot 10^{-5}$ 1/K. The body is initially at rest, it has an initial temperature $u_0=273.15$ K,  it is clamped at the bottom face, and it is heated at the top face with an ambient temperature $u_\infty=372.15$ K. The external heating introduced by the ambient temperature $u_\infty$ is modeled 
 as a prescribed flux, 
\begin{equation}
\label{eq:heat-conv-cond}
\qvec(u)\cdot\nvec = h(u-u_\infty),
\end{equation}
proportional to the difference between the current surface temperature $u$ and the ambient temperature $u_\infty$. The boundary condition \eqref{eq:heat-conv-cond} is the so-called heat convection boundary condition, and the proportionally constant $h$ is the so-called linear heat transfer coefficient. Its value is $h=\rho C \bar h$, with $\bar h=1\cdot 10^{-5}$ m/s.  
The other faces of the body are thermally insulated, i.e. homogeneous thermal Neumann conditions are applied. The problem is discretized in space with trilinear hexahedral finite elements. Four different meshes are considered in this example, see Table \ref{table:tsi-dani-meshes}. The time discretization is carried out with the one step-$\theta$ method with $\theta=2/3$ and time step $\Delta t=4\cdot 10^{-2}$ s.

\begin{table}[ht!]
\centering
\begin{tabular}{cccccc}
\toprule
Mesh Id& Processors & Structural field & Thermal field &  Total & Total/Processors   \\
\midrule
1&   2 &    63888 &   21296 &   85184 & 42592 \\
2&   8 &   235824 &   78608 &  314432 & 39304 \\
3&  32 &   998250 &  332750 & 1331000 & 41593 \\
4& 128 &  3951018 & 1317006 & 5268024 & 41156 \\
\bottomrule
\end{tabular}
\caption{Second Danilovskaya problem: Number of processors and number of degrees of freedom for the four meshes considered in the example.}
\label{table:tsi-dani-meshes}
\end{table}

The discrete non-linear problems to be solved at each time step are {\color{rev1} solved} with the monolithic Newton iteration \eqref{eq:TSI:Newton} and the resulting linear systems are solved using  GMRES preconditioned with the methods AMG(BGS), BGS({\color{rev1}AMG}) and SIMPLE(AMG) given in Table~\ref{table:preconds:tsi}.  For the outer Newton iterations, convergence is declared when the full residual fulfills $\Norm{\fmat^\mathrm{TSI}}_\mathrm{rms}<10^{-6}$ and the residual associated with the individual structural and thermal fields fulfill $\Norm{\fmat^\mathrm{S}}_\mathrm{rms}<10^{-8}$ and $\Norm{\fmat^\mathrm{Th}}_\mathrm{rms}<10^{-8}$, where $\Norm{\cdot}_\mathrm{rms}$ is the root mean square norm. The inner GMRES iteration is stopped when the residual of the linear problem is such that $\Norm{\rmat}_2/\Norm{\rmat^0}_2<10^{-8}$, where $\rmat^0$ stands for the residual associated with the initial guess given to the GMRES solver. Using these tolerances, three Newton iterations are required for achieving convergence at each time step.

The following parameters are considered for the studied preconditioners. The number of iterations for BGS and SIMPLE is taken as $k_\mathrm{max}=1$.  The underlying single field AMG methods are built using standard smoothed aggregation AMG \citep{VanekBrezinaMandel2001,VanekMandelBrezina1996} available in the MueLu \cite{ProkopenkoHuWiesnerSiefertTuminaro2014} library. The multigrid smoothers $\Svec^\ell_i$ for the single field problems are defined as an additive Schwartz domain decomposition using a damped Gauss-Seidel  iteration (damping $\omega=0.79$) as sub-domain solver. The coarse level solver is the so-called KLU direct solver given in the Amesos  package \citep{SalaStanleyHeroux2008}.

\begin{figure}[ht!]
\centering
\includegraphics[width=0.8\textwidth]{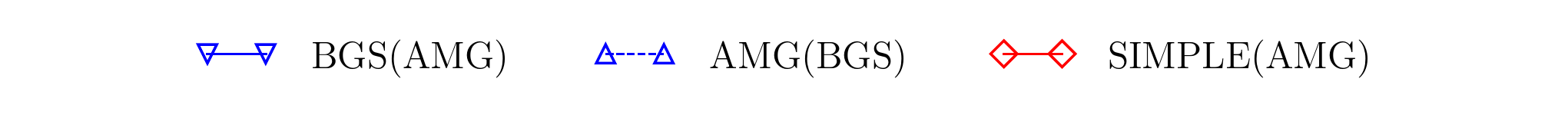}\\[0.2cm]

\begin{tabular}{cc}
\includegraphics[width=0.4\textwidth]{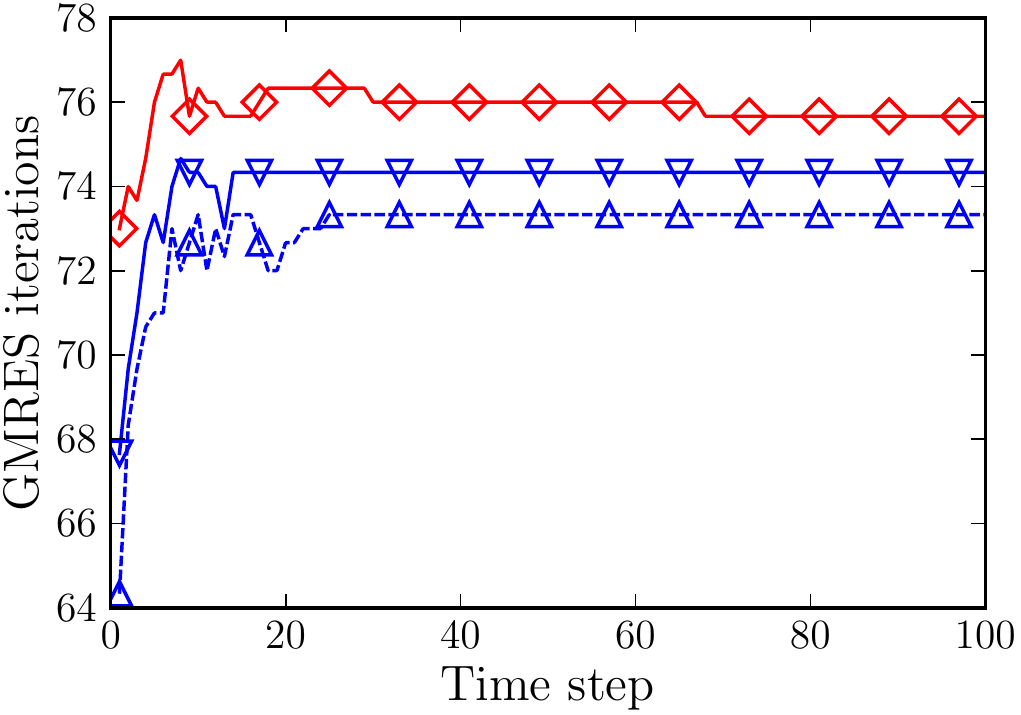} & \includegraphics[width=0.4\textwidth]{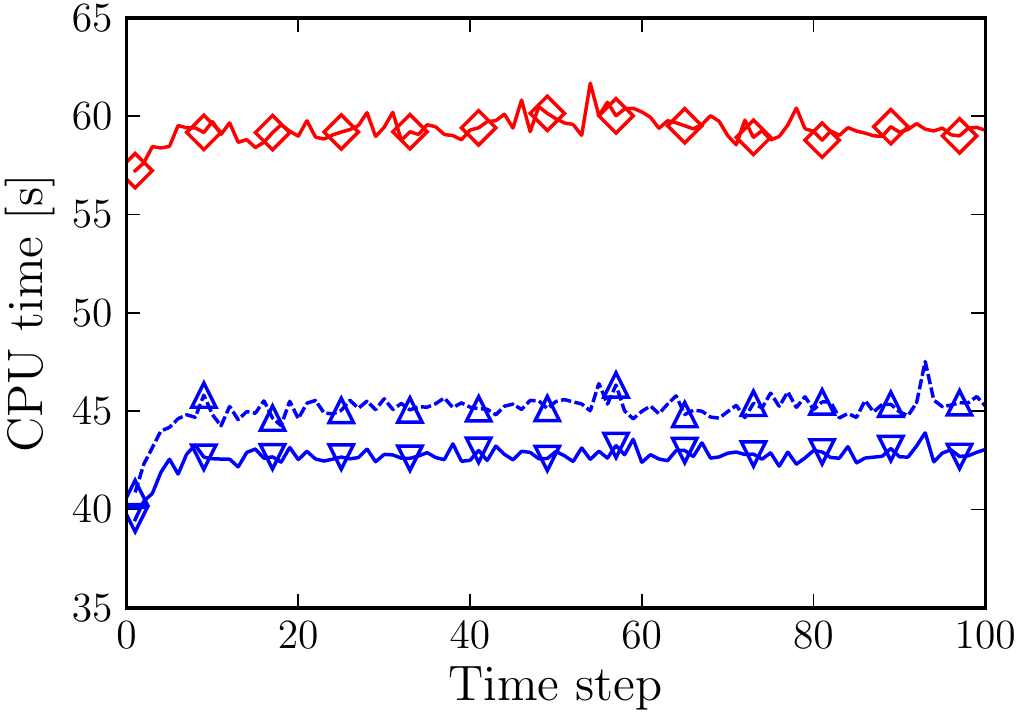}
\end{tabular}
\caption{Second Danilovskaya problem: Average of linear solver iterations per Newton step (left) and the corresponding CPU times (right). The results are for the finest discretization given in Table \ref{table:tsi-dani-meshes} (5.2 milion unknowns on 128 processors).}
\label{fig:tsi-dani-evol}
\end{figure}

\begin{figure}[ht!]
\centering
\includegraphics[width=0.8\textwidth]{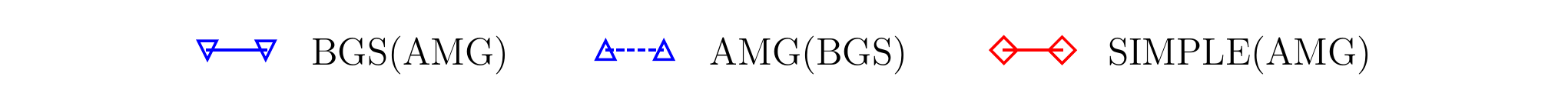}
\begin{tabular}{cc}
\includegraphics[width=0.4\textwidth]{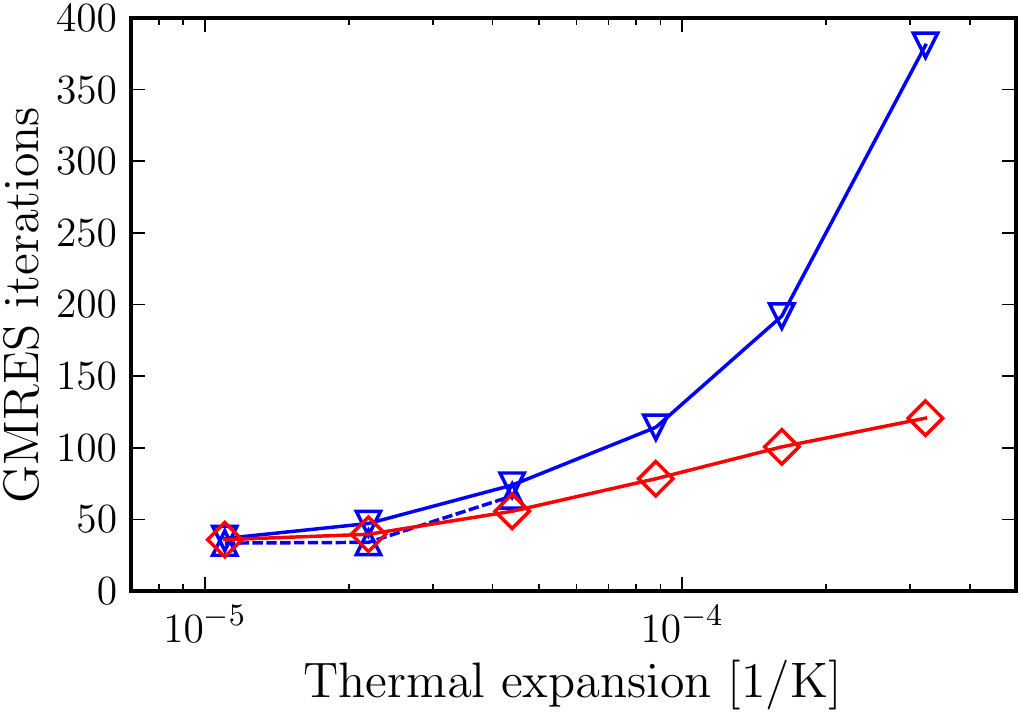} & \includegraphics[width=0.4\textwidth]{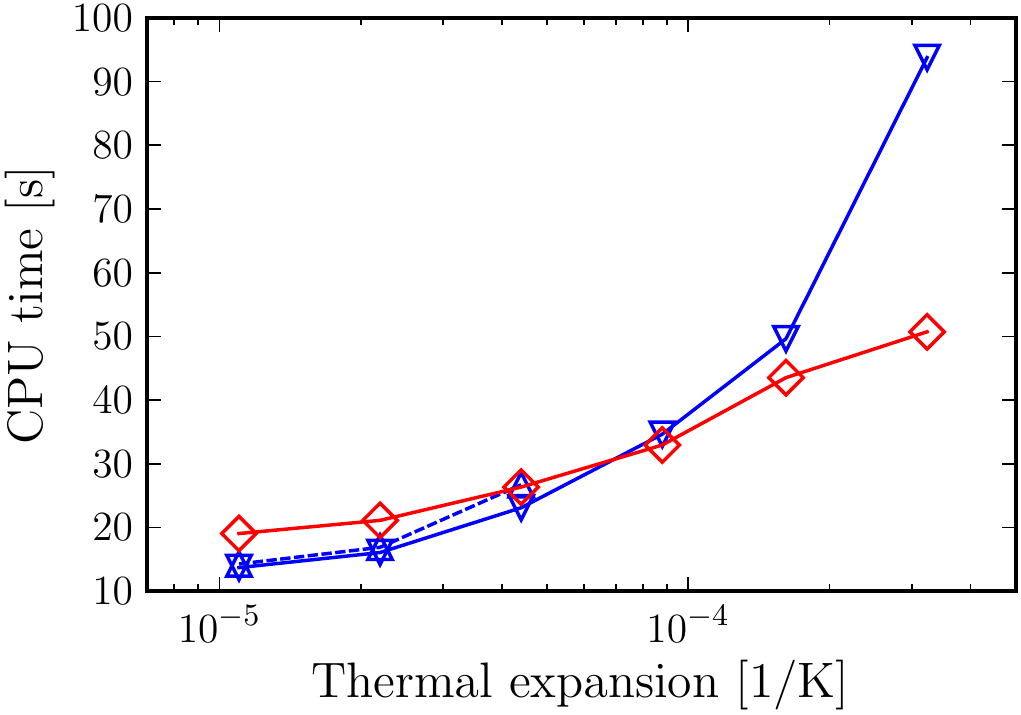}
\end{tabular}
\caption{Second Danilovskaya problem: Dependence of the performance of the    preconditioners  with respect to the amount of thermo-mechanical coupling. The results correspond to the first mesh given in Table \ref{table:tsi-dani-meshes}.}
\label{fig:tsi-dani-thexp}
\end{figure}

Figure \ref{fig:tsi-dani-evol} shows the average number of linear solver {\color{rev1}iterations} per Newton step and the associated total CPU time (including the setup of the preconditioners) for the three preconditioners considered here and for the finest mesh in Table \ref{table:tsi-dani-meshes}. The results show that the methods lead to a reasonable amount of GMRES iterations for this very fine mesh which means that the proposed methods are preconditioning the problem effectively. The best preconditioner in terms of iteration count is the AMG(BGS) method. This is what we expected since AMG(BGS) resolves the coupling at all multigrid levels. However, the BGS(AMG) method is {\color{rev1} clearly} better in terms of CPU times, which shows that the benefit of using AMG(BGS) (which is a more complex method) is not clearly justified for this example and the chosen parameters.
On the other hand, the SIMPLE(AMG) method is the worst in terms of iterations and CPU time. It might seem that SIMPLE(AMG) is not a suitable preconditioner for this type of problems. However, we show later in this example that for strong values of the thermal expansion coefficient $\alpha$, SIMPLE(AMG) is {\color{rev1} performing better than the} other two methods (see Figure \ref{fig:tsi-dani-thexp}).

\begin{figure}[ht!]
\centering
\includegraphics[width=0.8\textwidth]{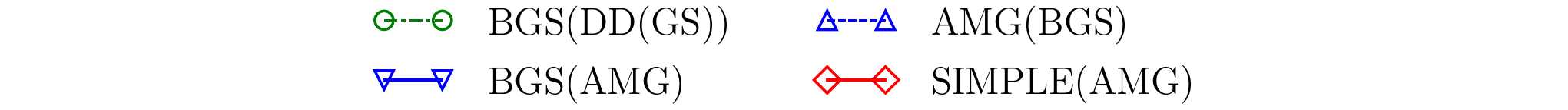}\\[0.2cm]

\begin{tabular}{cc}
\includegraphics[width=0.4\textwidth]{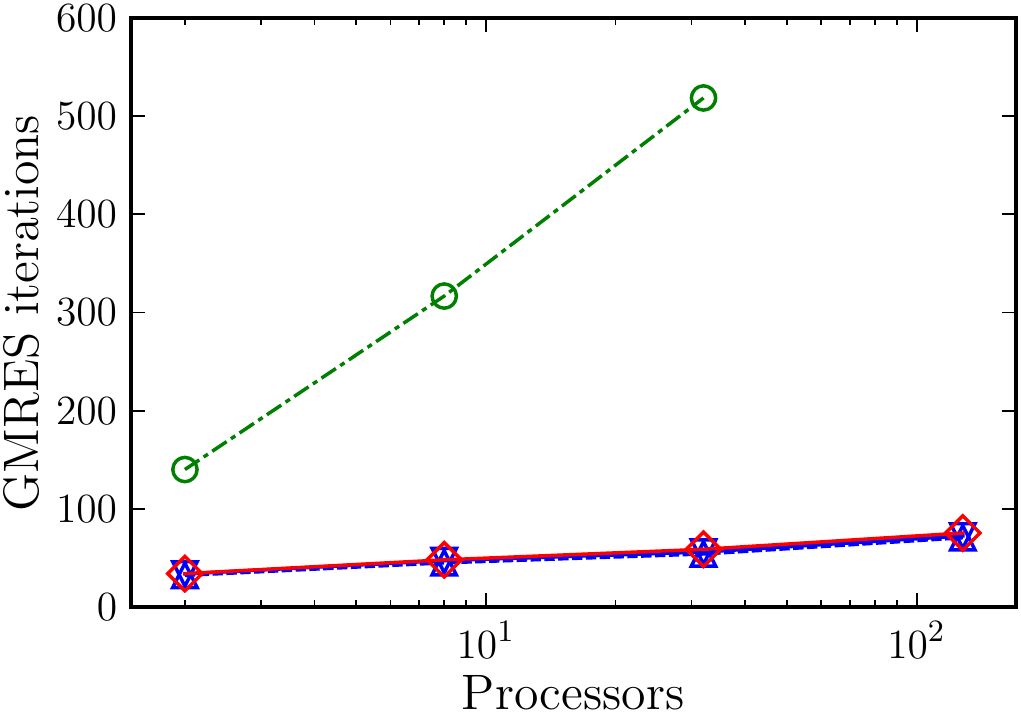} & \includegraphics[width=0.4\textwidth]{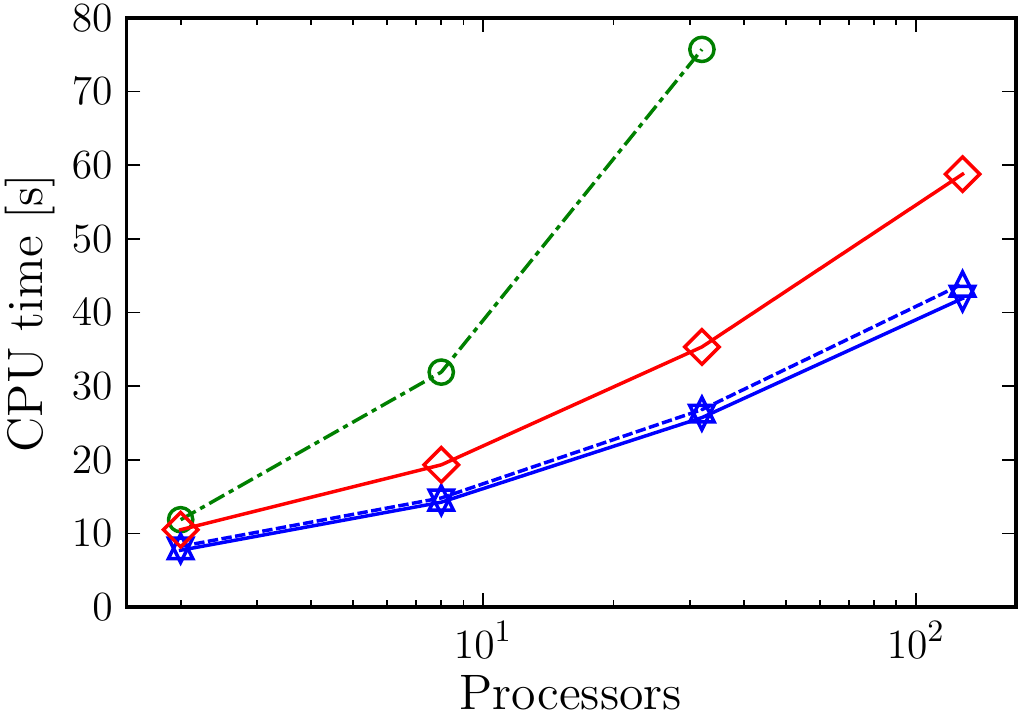}
\end{tabular}
\caption{Second Danilovskaya problem: Results of the weak scalability study.  The CPU times include the setup costs of the preconditioner. }
\label{fig:tsi-dani-weak-scal}
\end{figure}

Figure \ref{fig:tsi-dani-weak-scal}  and Tables  \ref{table:tsi-dani-weak-scal-iters} and \ref{table:tsi-dani-weak-scal-times}  report  a weak scalability test. The test is done using four meshes given in Table \ref{table:tsi-dani-meshes} and setting the number of processors such that the work per process is fairly the same in all cases. The proposed multigrid methods AMG(BGS), BGS(AMG) and SIMPLE(AMG) are  compared with a standard single-level preconditioner denoted here as BGS(DD(GS)). This single-level method is defined as an outer BGS iteration with standard additive Schwartz domain decomposition (DD)  methods for solving the resulting uncoupled problems. The local subdomain solvers in the domain decomposition are damped Gauss-Seidel  (GS) iterations with damping $\omega=0.79$. 
The single level method is one of the most simple parallel preconditioners that can be considered in such a problem, and therefore, it is considered here as a reference.

\begin{table}[ht!]
\centering
\begin{tabular}{cccccc}
\toprule
 Processors   & D.O.F.     &   BGS(DD(GS)) &  BGS(AMG)  &  AMG(BGS)  &  SIMPLE(AMG)\\
 \midrule
  2 &    85184  &        140  &        33  &            32  &            34   \\ 
  8 &   314432  &        316  &        47  &            45  &            48   \\  
 32 &  1331000  &        518  &        56  &            54  &            58   \\ 
128 &  5268024  &         --  &        72  &            70  &            75 \\  
\bottomrule
\end{tabular}
\caption{Second Danilovskaya problem: Average of linear solver iterations in the weak scalability test reported in Figure \ref{fig:tsi-dani-weak-scal}.}
\label{table:tsi-dani-weak-scal-iters}
\end{table}

\begin{table}[ht!]
\centering
\begin{tabular}{cccccc}
\toprule
 Processors   & D.O.F.     &   BGS(DD(GS)) &  BGS(AMG)  &  AMG(BGS)  &  SIMPLE(AMG)\\
 \midrule
  2 &    85184  &        11.8  &        7.7  &             8.2  &            10.4   \\ 
  8 &   314432  &        31.9  &       14.2  &            14.8  &            19.3   \\  
 32 &  1331000  &        75.7  &       25.7  &            26.8  &            35.3   \\ 
128 &  5268024  &         --  &        41.9  &            43.8  &            58.7 \\  
\bottomrule
\end{tabular}
\caption{Second Danilovskaya problem:  CPU times in the weak scalability test reported in Figure \ref{fig:tsi-dani-weak-scal}. }
\label{table:tsi-dani-weak-scal-times}
\end{table}

Figure \ref{fig:tsi-dani-weak-scal} and Tables \ref{table:tsi-dani-weak-scal-iters} and \ref{table:tsi-dani-weak-scal-times} show that the {\color{rev1} presented} multigrid methods AMG(BGS), BGS(AMG) and SIMPLE(AMG) lead to iteration counts and CPU times that mildly depend on the problem size. {\color{rev1} Obviously, this is not} the optimal multigrid behavior achieved in single-field elliptic problems (i.e. independence of the iteration count and CPU times with respect to the problem size), but it can be considered a very good scalability if one takes into account that the optimal multigrid {\color{rev1}behavior} is rarely reproduced in complex multiphysics. 
On the other hand, it is observed that the single-level method BGS(DD(GS)) is not scalable since the iteration count strongly depends on the problem size.  This is what we expected since {\color{rev1} it is well known that} single-level domain decomposition methods does not scale properly (cf. \cite{LinShadidTuminaroSalaHenniganPawlowski2010}). In conclusion, this scalability study shows that proposed preconditioners AMG(BGS), BGS(AMG) and SIMPLE(AMG)  have a decent scalability, and therefore, can be considered to solve this thermo-mechanical problem efficiently in parallel.

Finally, we investigate the performance of the preconditioners with respect to the thermal expansion coefficient $\alpha$, which controls the coupling between the structural and thermal fields (the coupling becomes stronger when $\alpha$ increases). 
As it is reported in Figure \ref{fig:tsi-dani-thexp}, the performance of the methods AMG(BGS) and BGS(AMG) strongly degrades when the thermal expansion coefficient $\alpha$ increases.
This is explained by noting that  the off-diagonal blocks of the system matrix in equation \eqref{eq:TSI:Monolithic} are proportional to  $\alpha$. Thus, when $\alpha$ grows, the off-diagonal blocks discarded by the BGS iteration become dominant and  affect negatively the performance of the method.  
The AMG(BGS) preconditioner is more sensitive to this effect than BGS(AMG). 
Note that for AMG(BGS), the convergence of the linear solver is achieved only up to the third value of the thermal expansion coefficient $\alpha=4.4\cdot 10^{-5}$ 1/K.
On the other hand, the SIMPLE(AMG) method is the most robust technique  with respect the thermal expansion coefficient $\alpha$. 
In conclusion, for mild values of the thermal expansion coefficient, {\color{rev1} the  BGS(AMG) method is preferred}, whereas for strong values, the best choice is  SIMPLE(AMG).

\subsection{TSI example: Rocket nozzle}

\begin{figure}[ht!]
\centering
\includegraphics[width=0.97\textwidth]{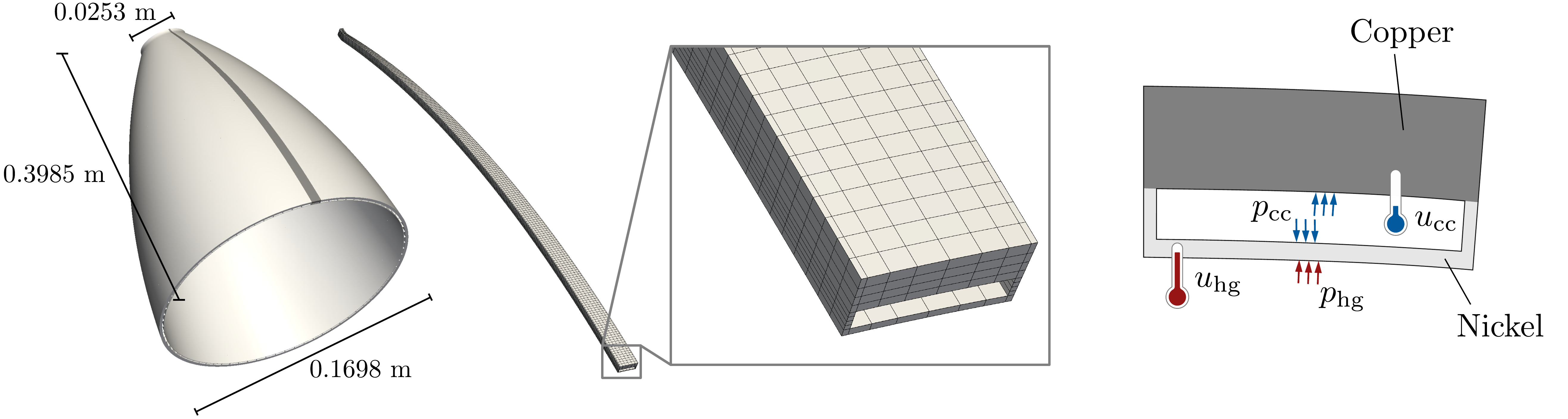}
\caption{Rocket nozzle example: Full geometry of the nozzle (left), computational domain including one cooling channel (center), and generic cross section of the computational domain with the applied thermo-mechanical loads (right). }
\label{fig:tsi-rocket-geom}
\end{figure}

The aim of this second example is to study the performance of the proposed preconditioners in a real-world problem. The example is taken form \citep{DanowskiGravemeierYoshiharaWall2013} and simulates the thermo-mechanical behavior of a rocked nozzle. 
The nozzle geometry  (see Figure \ref{fig:tsi-rocket-geom}) is inspired by the Vulcain rocket engine installed in the Ariane space launcher. The nozzle is $0.3985$ m in height and its cross section is approximated by a circular ring. The inner radius at the nozzle's inlet is $0.0253$ m and $0.1698$ m at the outlet.
The nozzle is equipped with 81 equidistant cooling channels which traverse the body longitudinally. Thanks to the symmetry of the problem, only a 1/81-th of the  nozzle is meshed and simulated by introducing the required symmetry conditions. The computational domain is a sector of the nozzle of angle $360^\circ/81$ which includes only one cooling channel, see Figure  \ref{fig:tsi-rocket-geom} (center).

The rocket nozzle is made of two materials: an external copper alloy body and an internal nickel jacket, see Figure \ref{fig:tsi-rocket-geom} (right). The copper alloy has  Young's modulus  $E=148$ GPa, Poisson's ratio $\nu=0.3$,  density $\rho=9.13\cdot 10^3$ kg/m$^3$, heat conductivity $k=310$ W/(m $\cdot$ K), heat capacity  $C=373$ J/(kg $\cdot$ K) and coefficient of thermal expansion $\alpha=1.72\cdot 10^{-5}$ 1/K. On the other hand, the nickel jacked has the following properties: $E=193$ GPa, $\nu=0.3$, $\rho=8910$ kg/m$^3$, $k=75$ W/(m $\cdot$ K),  $C=444$ J/(kg $\cdot$ K) and    $\alpha=1.22\cdot 10^{-5}$ 1/K.

\begin{figure}[ht!]
\centering
\begin{tabular}{cc}
\includegraphics[width=0.4\textwidth]{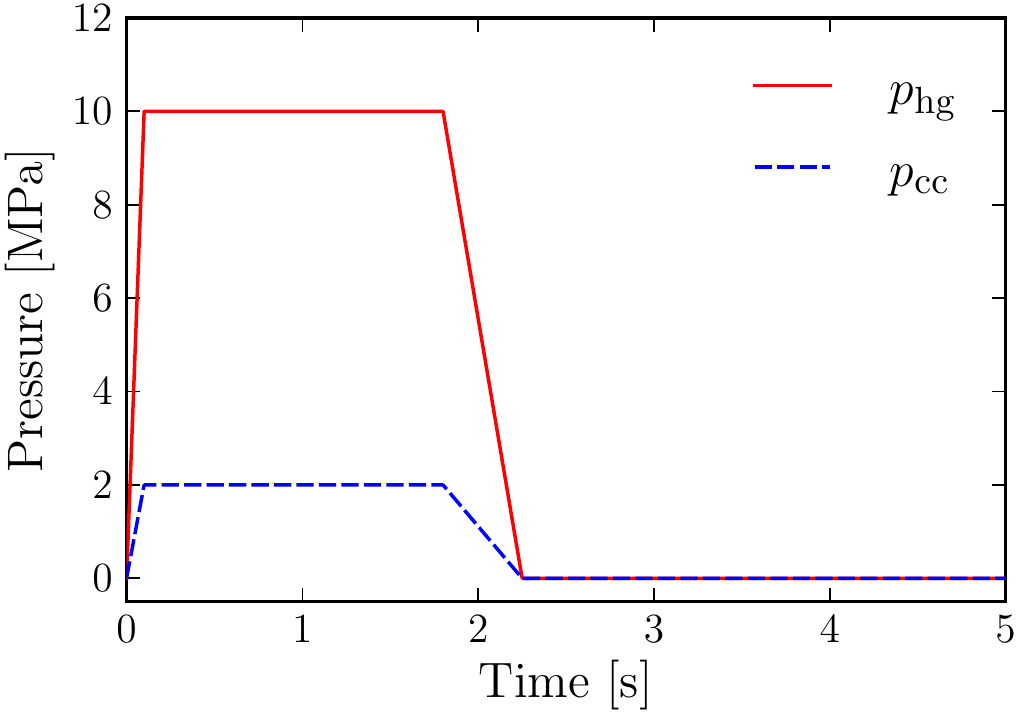} & \includegraphics[width=0.4\textwidth]{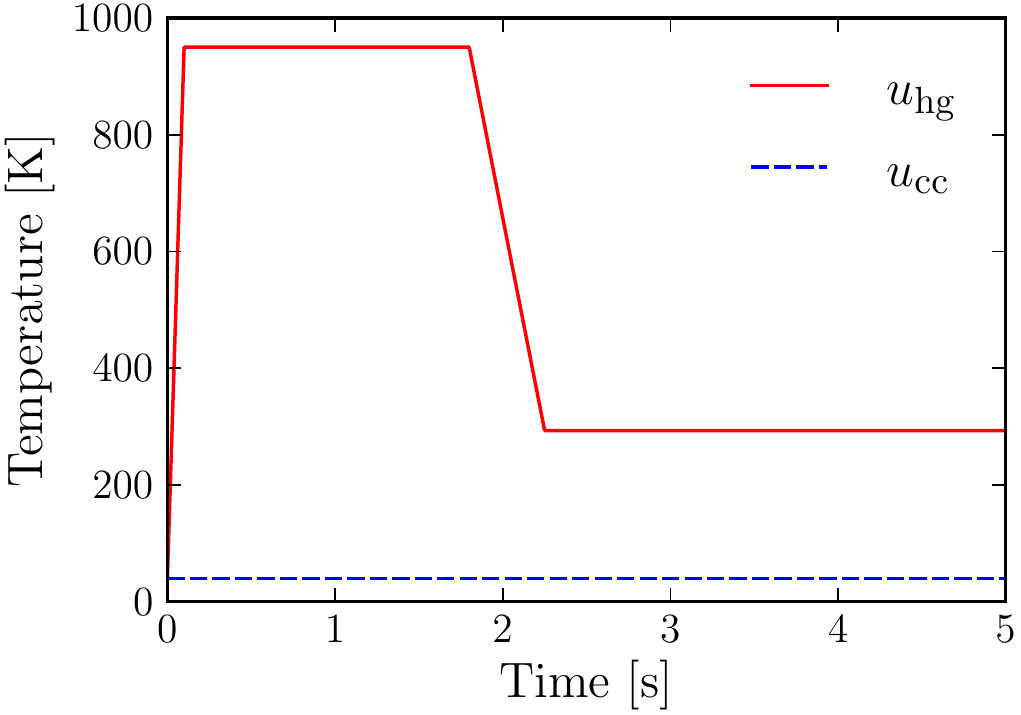}
\end{tabular}
\caption{Rocket nozzle example: Temporal evolution of the thermo-mechanical loads (pressure and temperature) associated with the exhausted hot gas ($p_\mathrm{hg}$, $u_\mathrm{hg}$) and with the liquid filling the cooling channels ($p_\mathrm{cc}$, $u_\mathrm{cc}$).}
\label{fig:tsi-nozzle-loads}
\end{figure}

\begin{figure}[ht!]
\centering
\includegraphics[width=0.93\textwidth]{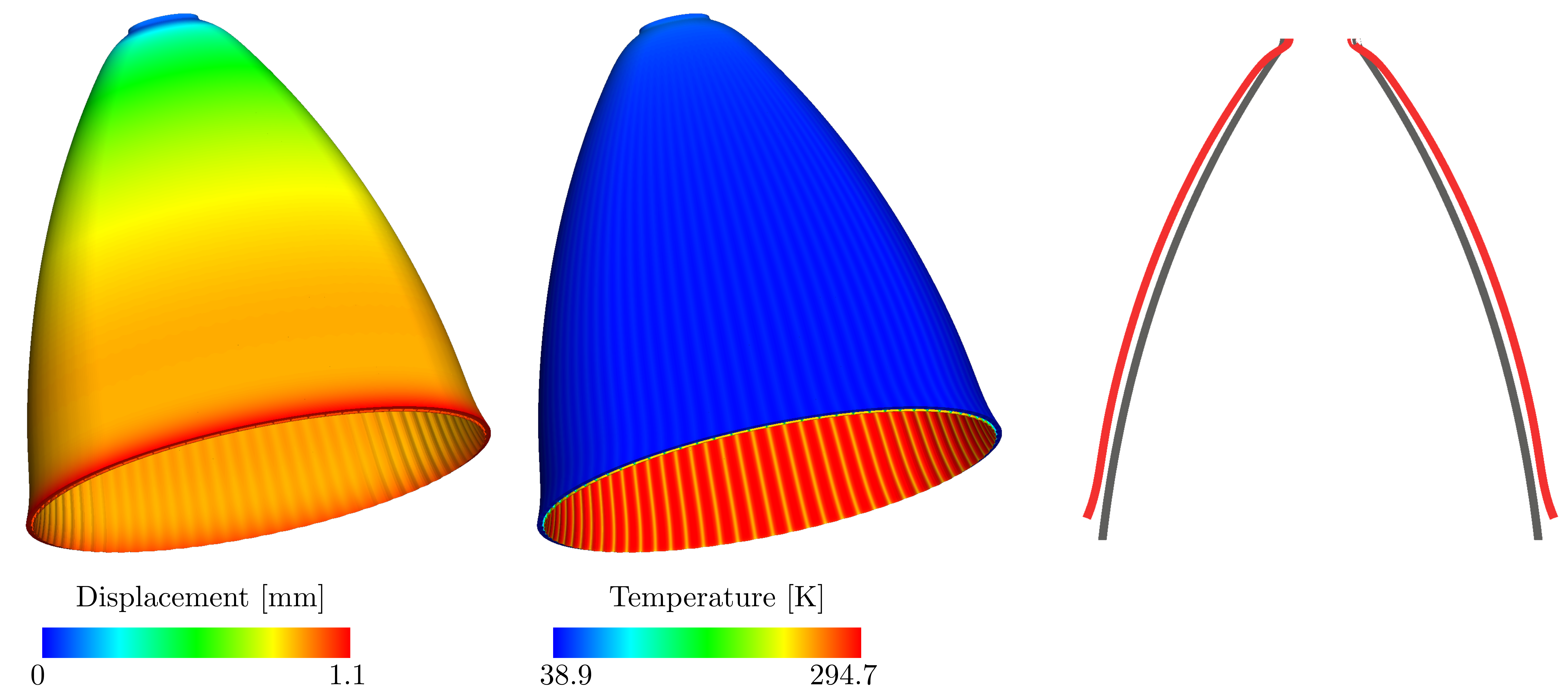}
\caption{Rocket nozzle example: Deformation of the nozzle (left), temperature distribution (center) and original and deformed longitudinal section (right). The results are given at time $t=1$ s and the deformation is magnified 20 times. }
\label{fig:tsi-rocket-sol}
\end{figure}

The rocked nozzle is initially at rest and cooled down with {\color{rev1} a} uniform initial temperature of $40$ K. The {\color{rev1}thermo}-mechanical loads acting on the body are the high pressure and temperature of the exhausted hot gases (namely $p_\mathrm{hg}$ and $u_\mathrm{hg}$) and the pressure and low temperature of the fluid filling the cooling channels (namely $p_\mathrm{cc}$ and $u_\mathrm{cc}$). The  temporal evolution of  $p_\mathrm{hg}$, $u_\mathrm{hg}$, $p_\mathrm{cc}$ and $u_\mathrm{cc}$ is  taken from \citep{AryaArnold1992} corresponding to a prototypical  loading cycle for rocked nozzles, see Figure \ref{fig:tsi-nozzle-loads}. As in the previous example, the effect of the external temperatures $u_\mathrm{hg}$ and $u_\mathrm{cc}$ is modeled with the heat convection boundary condition \eqref{eq:heat-conv-cond}. 
The linear heat transfer coefficient $h$ is $h_\mathrm{hg}=32$ kW/(m$^2$K) for the hot gas, and   $h_\mathrm{cc}=100$ kW/(m$^2$K) for the cooling channel.
Additionally, the inlet of the nozzle is fixed in longitudinal direction in order to take into account the combustion chamber which is attached to it.
The total simulation time is $t=5$ s. A screen-shot of the deformation of the nozzle and the temperature field is given in Figure  \ref{fig:tsi-rocket-sol}.

The performance of the preconditioners is studied with a weak scalability test. To this end, several finite element meshes are considered and the number of processors is chosen {\color{rev1} such} that the work per processor is similar for all the meshes, see Table  \ref{table:tsi-rocket-discr}.
As in the previous example, tri-linear hexahedral finite elements are used for the space discretization and the  one step-$\theta$ with $\theta=2/3$ for the time discretization.  The time step is chosen as $\Delta t=5 \cdot 10^{-2}$ s.  The Newton iterations are stopped when the full non-linear residuum is $\Norm{\fmat^\mathrm{TSI}}_\mathrm{rms}<10^{-6}$, and the entries of the non-linear residuum corresponding to the structural and thermal fields are  $\Norm{\fmat_\mathrm{S}}_\mathrm{rms}<10^{-6}$ and $\Norm{\fmat_\mathrm{T}}_\mathrm{rms}<10^{-6}$ respectively. The GMRES iteration is stopped when  the linear residuum fulfills $\Norm{\rmat}_2/\Norm{\rmat^0}_2<10^{-7}$. Using these tolerances, between 2 and 3 Newton iterations are required for achieving convergence at each time step. The preconditioners AMG(BGS), BGS(AMG), SIMPLE(AMG) and BGS(DD(GS)) are the same as in the previous example.

\begin{table}[ht!]
\centering
\begin{tabular}{cccccc}
\toprule
Mesh Id& Processors & Solid & Thermal field &  Total & Total/Processors   \\
\midrule
1&    4 &    47025 &  15675 &   62700 &  15675 \\
2&   32 &   307046 & 102582 &  410328 &  12822 \\
3&  256 &  2189124 & 729708 & 2918832 & 11401 \\
\bottomrule
\end{tabular}
\caption{Rocket nozzle example: Number of processors and number of degrees of freedom for the three meshes considered in the example.}
\label{table:tsi-rocket-discr}
\end{table}

\begin{figure}[ht!]
\centering
\includegraphics[width=0.8\textwidth]{figures/legend.pdf}\\[0.2cm]

\begin{tabular}{cc}
\includegraphics[width=0.4\textwidth]{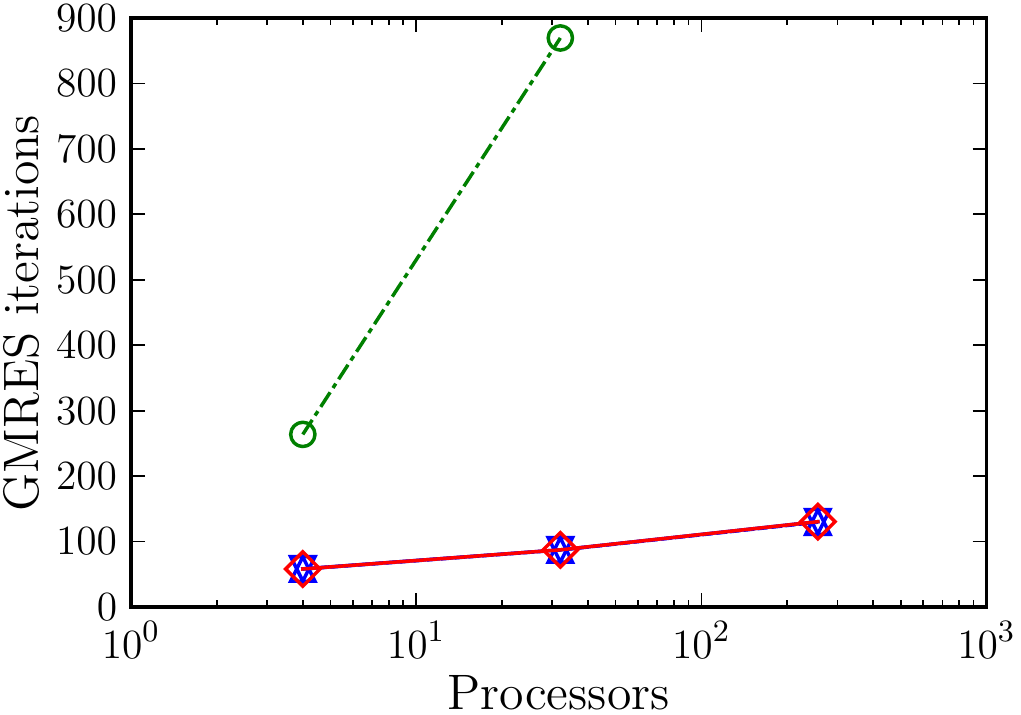} & \includegraphics[width=0.4\textwidth]{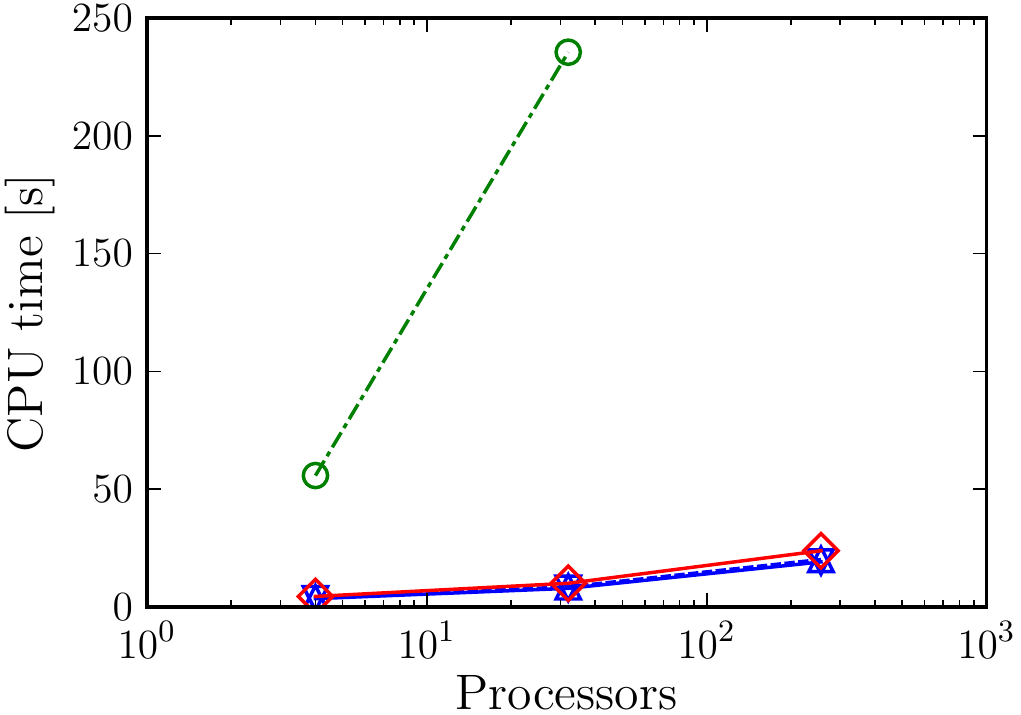}
\end{tabular}
\caption{Rocket nozzle example: Results of the weak scalability study.  The CPU times include the setup costs of the preconditioner.}
\label{fig:tsi-rocket-weak-scal}
\end{figure}

\begin{table}[ht!]
\centering
\begin{tabular}{cccccc}
\toprule
 Processors   & D.O.F.     &   BGS(DD(GS)) &  BGS(AMG)  &  AMG(BGS)  &  SIMPLE(AMG)\\
 \midrule
   4 &    62700  &        263  &         58  &            57  &          58     \\ 
  32 &   410328  &        869  &         87  &            87  &          87     \\  
 256 &  2918832  &         --  &        130  &            129  &          130     \\ 
\bottomrule
\end{tabular}
\caption{Rocket nozzle example: Average of linear solver iterations in the weak scalability test reported in Figure~\ref{fig:tsi-rocket-weak-scal}.}
\label{table:tsi-rocket-weak-scal-iters}
\end{table}

\begin{table}[ht!]
\centering
\begin{tabular}{cccccc}
\toprule
 Processors   & D.O.F.     &   BGS(DD(GS)) &  BGS(AMG)  &  AMG(BGS)  &  SIMPLE(AMG)\\
 \midrule
   4 &    62700  &         55.7  &         3.5  &            3.9  &          4.5     \\ 
  32 &   410328  &        235.4  &         7.9  &            8.5  &          10.1     \\  
 256 &  2918832  &           --  &        19.1  &           20.1  &          23.9     \\ 
\bottomrule
\end{tabular}
\caption{Rocket nozzle example:  CPU times in the weak scalability test reported in Figure \ref{fig:tsi-rocket-weak-scal}.}
\label{table:tsi-rocket-weak-scal-time}
\end{table}

The results of the weak scalability study are reported in Figure \ref{fig:tsi-rocket-weak-scal} and in Tables \ref{table:tsi-rocket-weak-scal-iters} and \ref{table:tsi-rocket-weak-scal-time}. The performance of the preconditioners in this complex three dimensional example is similar to the behavior observed in the previous {\color{rev1}benchmark} test (Section \ref{sec:TSI:example:danilovsakaya}).  Again, the multigrid methods AMG(BGS), BGS(AMG) and SIMPLE(AMG) have a very good scalability as the iteration count and CPU time of the linear solver increase only mildly with the problem size. On the other hand, the single-level method BGS(DD(GS)) based on a standard additive Schwarz domain decomposition is not a scalable method as previously observed.
Note that the performance of the single-level method BGS(DD(GS)) is particularly poor in this complex example which demonstrates that the multigrid preconditioners  AMG(BGS), BGS(AMG) and SIMPLE(AMG) are doing a particularly efficient job in this challenging setting. In conclusion, the  {\color{rev1} multigrid} methods  AMG(BGS), BGS(AMG) and SIMPLE(AMG) have a good scalability also for this challenging example, which shows that the methods are able to solve this complex real-world problem efficiently.

\subsection{FSI example: Pressure wave}

\begin{figure}[ht!]
\centering
\includegraphics[width=0.9\textwidth]{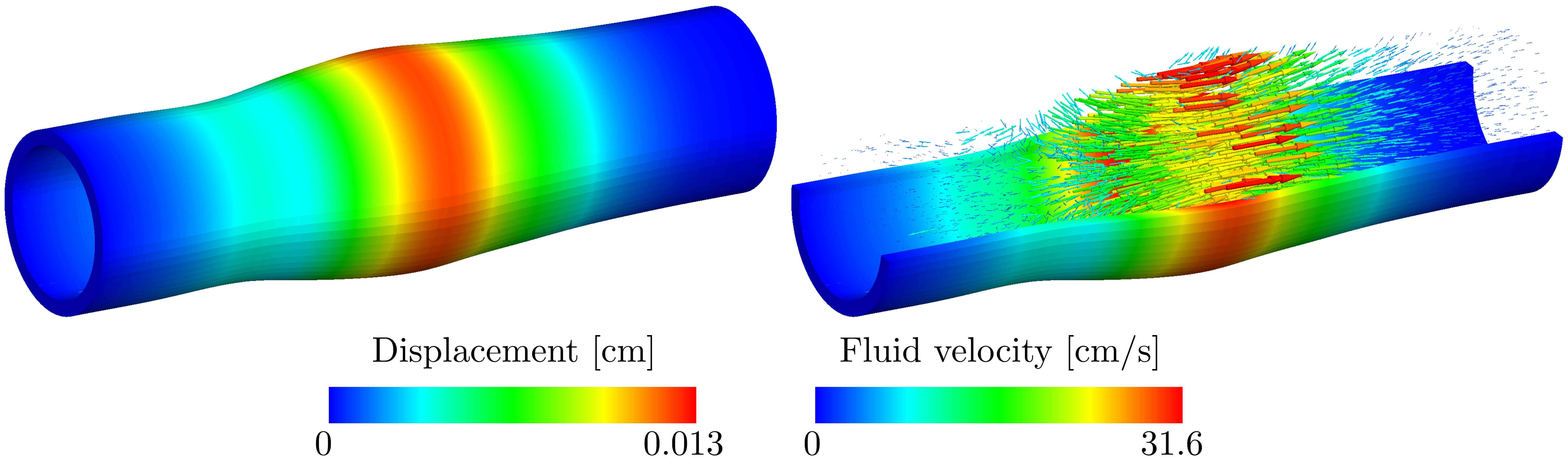}
\caption{FSI example: Pressure wave traveling through a flexible tube. 
The deformation is magnifyied 10 times.   }
\label{fig:fsi:cyl:solution}
\end{figure}

The aim of this example is to compare the performance of the proposed general purpose preconditioners with respect to the FSI-specific preconditioners proposed by \citet{GeeKuettlerWall2011}.
The study is carried out for the well known benchmark test introduced by \citet{GerbeauVidrascu2003} consisting in a fluid wave traveling through a flexible tube which mimics hemodynamic conditions, see Figure \ref{fig:fsi:cyl:solution}. The tube has a length of $5\cdot 10^{-2}$ m, an inner radius of $5\cdot 10^{-3}$ m and an outer radius of $6\cdot 10^{-3}$ m. The structural model used for the tube is a geometrically non-linear Saint Venant-Kirchhoff material with Young's modulus  $E=3\cdot 10^5$ Pa, Poisson's ratio $\nu=0.3$ and density $\rho^\mathrm{s}=1.2\cdot 10^3$ Kg/m$^3$. The tube is filled with a Newtonian fluid with density $\rho^\mathrm{f}=1\cdot 10^3$ Kg/m$^3$ and viscosity $\mu = 3\cdot 10^{-3}$ Pa $\cdot$ s. Both the tube and the contained fluid are at rest at the initial time. The tube is clamped at both ends and the fluid is initially loaded with a surface traction of  $1.33\cdot 10^3$ Pa on one side for a time interval of $3\cdot 10^{-3}$ s. This action results in a pressure wave traveling through the tube, see Figure  \ref{fig:fsi:cyl:solution}.

The problem is discretized in space using four different finite element meshes, see Table \ref{table:discr:fsi:cyl}. Trilinear hexahedral elements are used in the structural domain and the fluid is discretized with stabilized trilinear hexahedral elements using the same interpolation for both velocity and pressure. The time discretization is performed using the trapezoidal rule for the structural domain and with the one step-$\theta$ method for the fluid ($\theta=2/3$). The time step is $\Delta t=10^{-4}$ s.

\begin{table}[ht!]
\centering
\begin{tabular}{ccccccc}
\toprule
Mesh Id& Processors & Solid & Fluid & Grid & Total & Total/Processors   \\
\midrule
1& 1  &   7560  &   23660 &  13965 &   45185 & 45185\\
2& 4  &  24708  &   97128 &  60492 &  182328 & 45582\\
3& 16 &  95580  &  372880 & 247800 &  716260 & 44766\\
4& 64 & 396000  & 1434752 & 996864 & 2827616 & 44181\\
\bottomrule
\end{tabular}
\caption{FSI example: Number of processors and number of degrees of freedom for the four meshes considered in the example. }
\label{table:discr:fsi:cyl}
\end{table}

The non-linear FSI problems associated with each time step are solved with the monolithic Newton iteration \eqref{eq:FSI:Newton} and using a preconditioned GMRES method for solving the resulting linear systems. 
  The  outer Newton iteration is stopped when the full FSI-residual fulfills $\Norm{\fmat^\mathrm{FSI}}_\mathrm{rms}<10^{-6}$ and $\Norm{\fmat^\mathrm{FSI}}_\infty<10^{-6}$  and the linear solver is stopped when the linear residual fulfills $\Norm{\rmat}_2/\Norm{\rmat^0}_2<10^{-10}$. 

The preconditioners studied here (Table \ref{table:preconds:fsi}) require building AMG methods for each  of the FSI fields. For the structural field and the ALE field, the underlying multigrid hierarchies are built using standard smoothed aggregation AMG \citep{VanekBrezinaMandel2001,VanekMandelBrezina1996}. For the fluid field,  a Petrov-Galerkin \citep{SalaTuminaro2008} AMG method is considered which is specially designed for the non-symmetric structure of the fluid equations. The multigrid smoothers are chosen as an additive Schwartz domain decomposition  using a damped Gauss-Seidel  iteration as subdomain solver, see Table \ref{table:fsi-cyl-amg}. This smoother is referred to as  DD(GS($\omega$)), where $\omega$ is the damping factor of the Gauss-Seidel iteration. The coarse level solver is the KLU direct solver given in the Amesos  package \citep{SalaStanleyHeroux2008}. 

\begin{table}[ht!]
\centering
\begin{tabular}{ccccc}
\toprule

Field & Smoothers & Coarsest level solver & Transfer operators & \# of levels\\
\midrule
Solid & DD(GS($0.79$))& KLU & Smoothed aggregation AMG   & 3\\
Fluid & DD(GS($0.69$))& KLU & Petrov-Galerkin AMG        & 3\\
ALE   & DD(GS($0.79$))& KLU & Smoothed aggregation AMG   & 3 \\
\bottomrule
\end{tabular}
\caption{FSI example: Details of the AMG methods used to build the preconditioners. }
\label{table:fsi-cyl-amg}
\end{table}

\begin{figure}[ht!]
\centering
\includegraphics[width=0.8\textwidth]{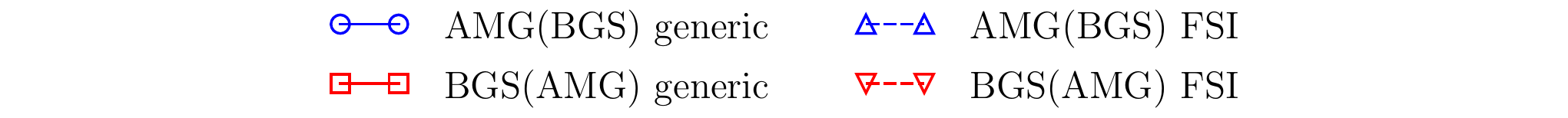}\\[0.2cm]

\begin{tabular}{cc}
\includegraphics[width=0.4\textwidth]{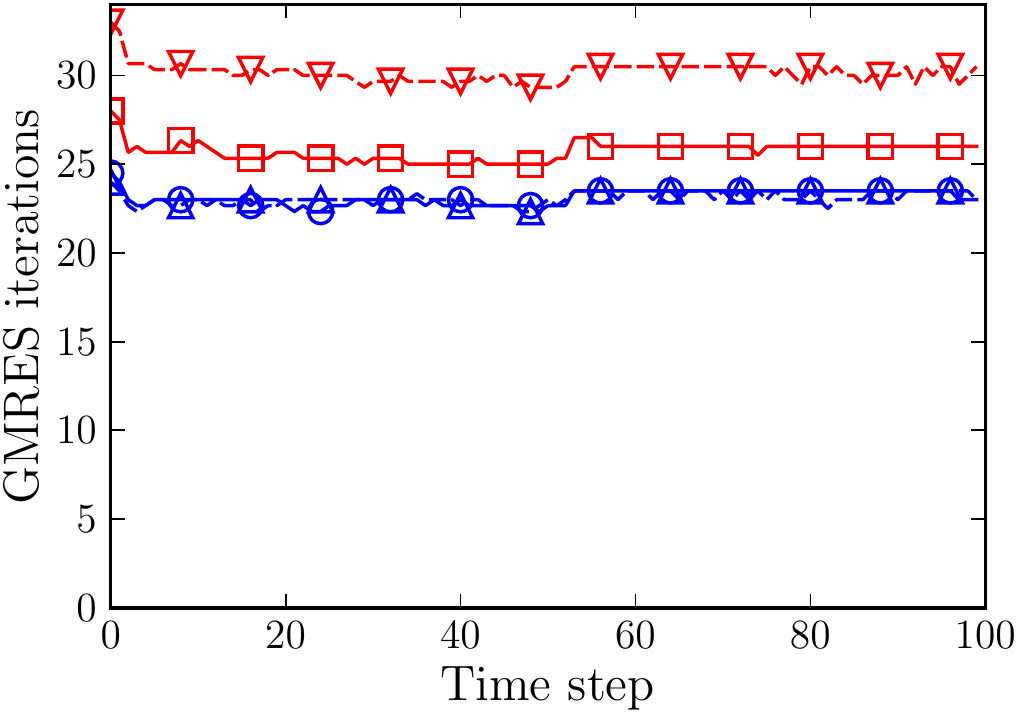} & \includegraphics[width=0.4\textwidth]{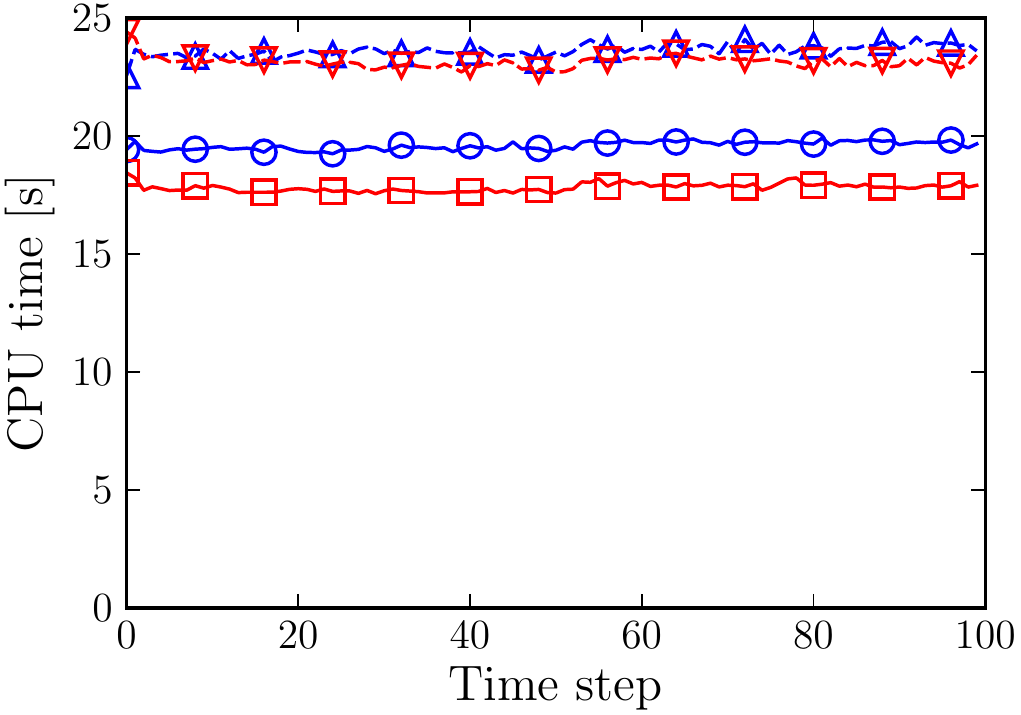}
\end{tabular}
\caption{FSI example: Performance of the generic and FSI-specific preconditioners. The CPU times include the setup cost of the preconditioner. The results are for the finest discretization given in Table \ref{table:discr:fsi:cyl} (2.8 {\color{rev1}million} unknowns on 64 processors).  }
\label{fig:fsi-cyl-evol}
\end{figure}

Figure \ref{fig:fsi-cyl-evol} shows the  performance of the general purpose preconditioners BGS(AMG) and AMG(BGS) and their corresponding FSI-specific implementations proposed by \citet{GeeKuettlerWall2011}. Note that both implementations lead to similar CPU time which show that using a generic implementation (which is not optimized for this particular problem) is not affecting the performance of the solver. In fact, the generic implementation has slightly better CPU times, which can be  due to the fact that different packages are considered for building the underlying AMG {\color{rev1}methods}:
the FSI-specific implementation  uses ML \citep{GeeSiefertHuTuminaroSala2006}  and the general implementation proposed here uses MueLu \citep{ProkopenkoHuWiesnerSiefertTuminaro2014} which is a more recent multigrid package. This, shows that the performance of the methods is driven by the chosen AMG libraries (which are black boxes here), and not by the implementation of the coupling iteration.

\begin{figure}[ht!]
\centering
\includegraphics[width=0.8\textwidth]{figures/fsi_cyl_evol_legend.pdf}\\[0.2cm]

\begin{tabular}{cc}
\includegraphics[width=0.4\textwidth]{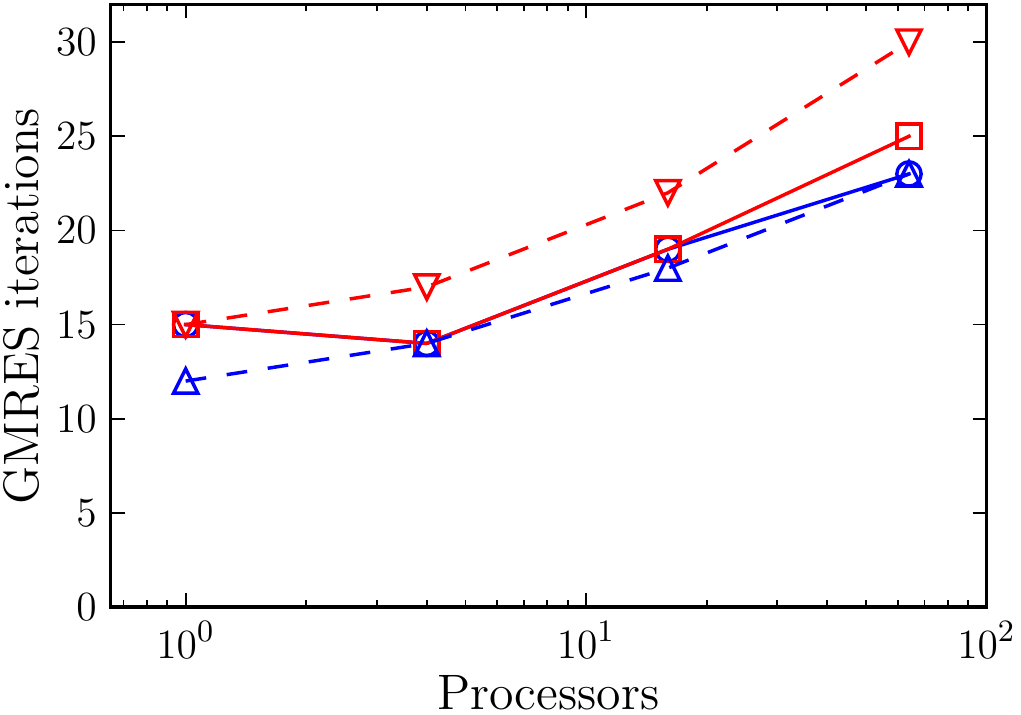} & \includegraphics[width=0.4\textwidth]{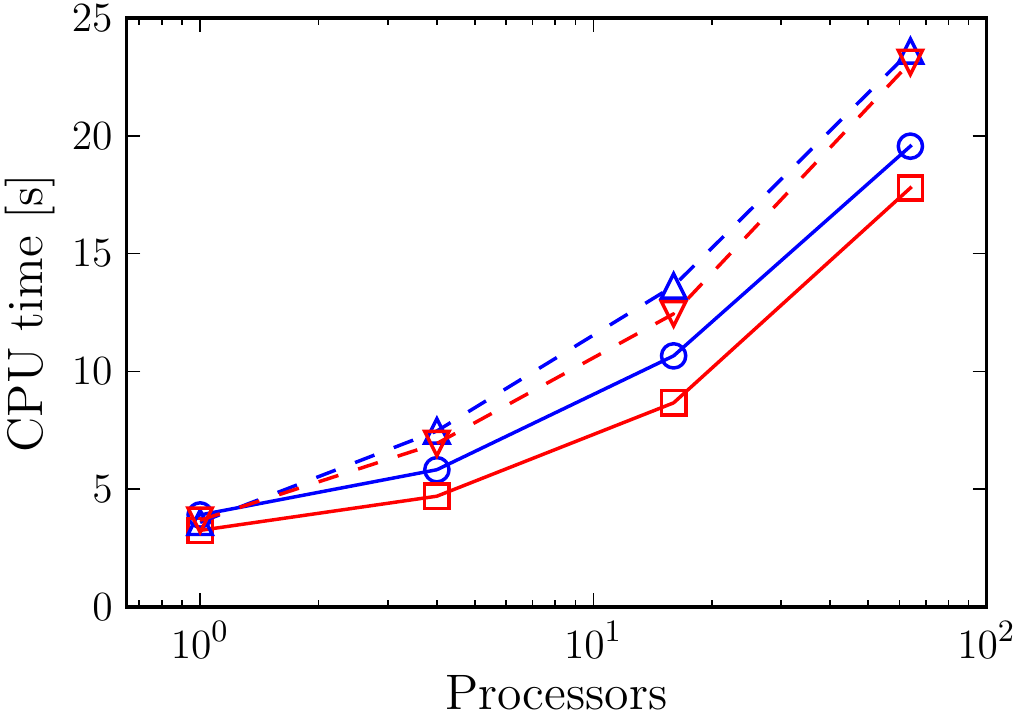}
\end{tabular}
\caption{FSI example: Results of the weak scalability study.  The CPU times include the setup costs of the preconditioner. }
\label{fig:fsi-cyl-weack-scal}
\end{figure}

The same comparison is carried out using a weak scalability study. To this end, four meshes are considered (see Table \ref{table:discr:fsi:cyl}) and the number of processors are chosen such that the amount of work per processor is  fairly constant (circa 45000 unknowns per processor). The results (Figure~\ref{fig:fsi-cyl-weack-scal}) show that the studied methods have a good  scalability also for this FSI example. Again, the observed scalabilities do not correspond to the optimal multigrid behavior expected in elliptic problems, but they can be considered decent results for this complex application. Note also that both the generic and the FSI-specific implementations lead to similar  CPU times, which confirms that the general purpose preconditioners  are efficient methods.

\subsection{Patient-specific lung example}

\begin{figure}[ht!]
\centering

\begin{tabular}{cc}
\includegraphics[width=0.31\textwidth]{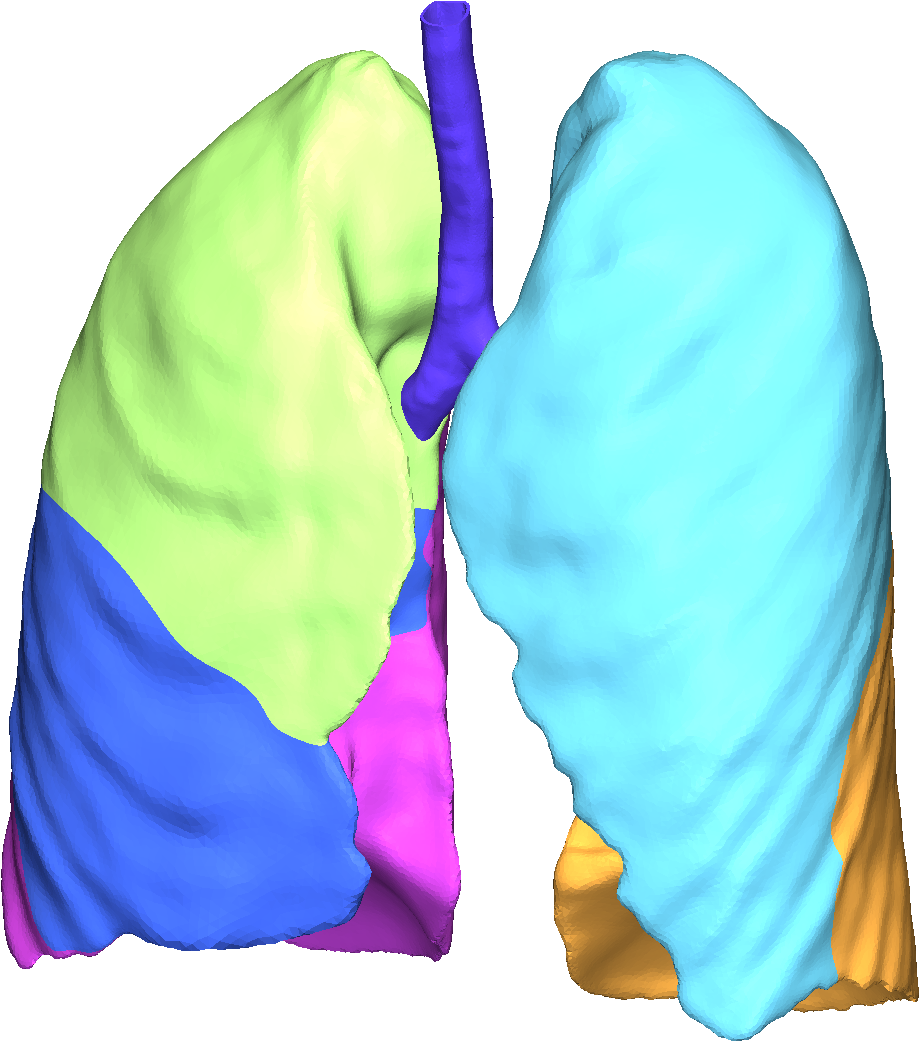} & \includegraphics[width=0.31\textwidth]{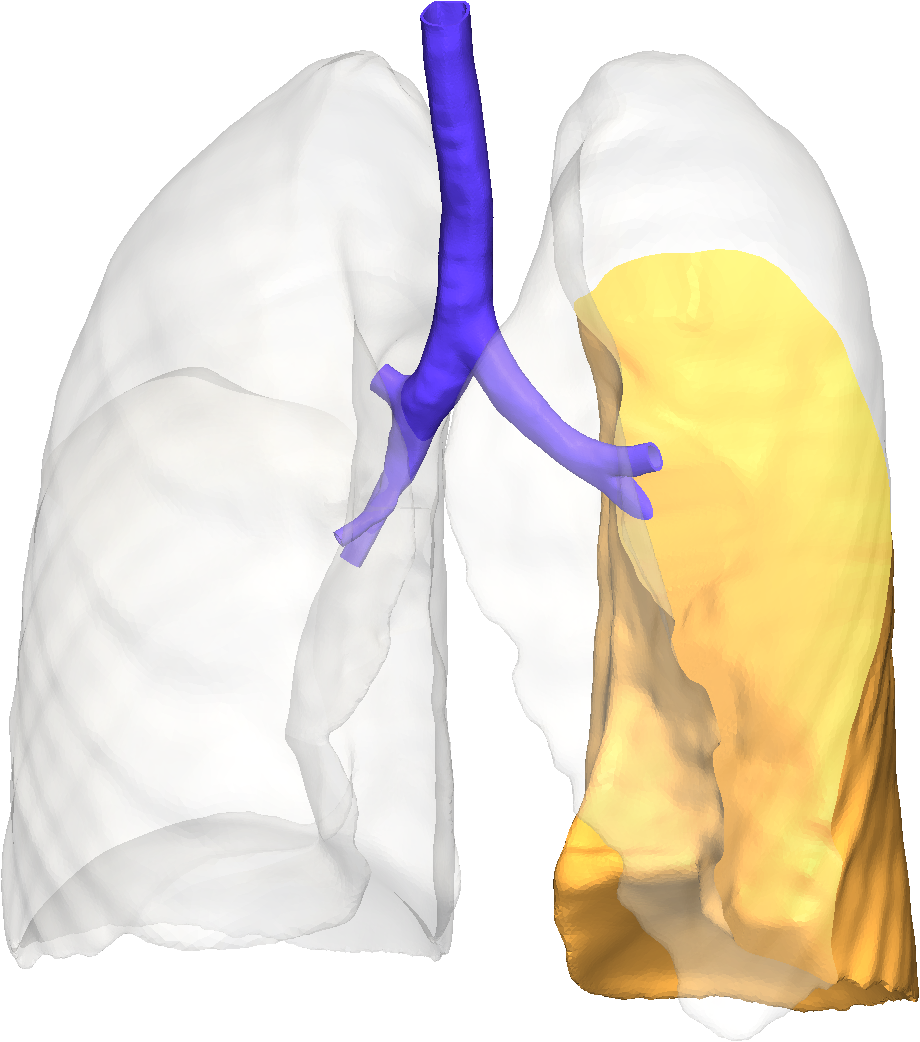}
\end{tabular}

\caption{Patient-specific lung example: Lung geometry with the airway tree reconstructed up-to the 3rd generation. The left picture shows the different regions in which the lung is sub-divided according to the model by \citet{YoshiharaIsmailWall2013,YoshiharaRothWall2015}. The right picture shows the airway tree and one of the these regions.}
\label{fig:lung-geom}
\end{figure}

Finally, we study the applicability of the proposed  methods to the lung model proposed by \citet{YoshiharaIsmailWall2013,YoshiharaRothWall2015} (see Section \ref{sec:lung:theory}) using a real-world patient-specific example.
The problem geometry corresponds to the airway tree and the lung lobes (see Figure \ref{fig:lung-geom}) of a healthy male subject, which is obtained from computer tomography (CT) images provided by the Medical Office ``Diagnostische Radiologie'' in Stuttgart, Germany. The geometry of the air tree is reconstructed up to the 3rd generation leading to 5 artificial outlets that have to be coupled to the lung tissue. This is achieved by subdividing the lung lobes into several sub-domains associated with each one of the ending branches (see Figure \ref{fig:lung-geom}). These regions represent the volumes enclosed by the boundaries $\Gamma_\bullet$  in Section  \ref{sec:lung:theory}.

\begin{figure}[ht!]
\centering
\includegraphics[width=0.94\textwidth]{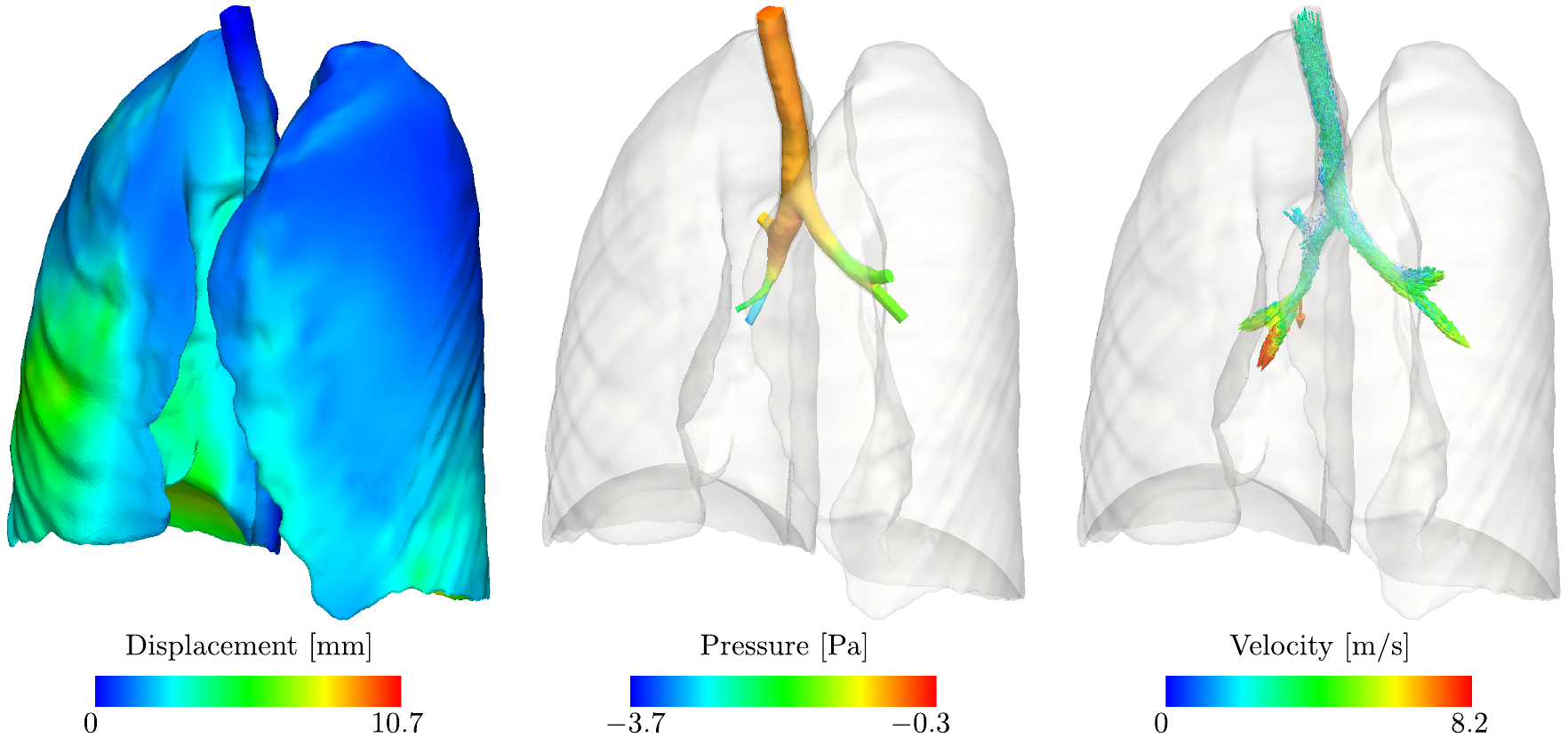}
\caption{Example 2: Numerical solution of the lung model consisting in the structural displacement (left), fluid pressure (center) and fluid velocity (right) at time $t=0.75$ s.}
\label{fig:lung-sol}
\end{figure}

The lung tissue is modeled with a Neo-Hookean solid with young's modulus $E=6.75$ kPa, Poisson's ration $\nu=0.3$ and density $\rho^\mathrm{s}=1.2\cdot 10^3$ Kg/m$^3$ which are determined experimentally in \cite{RauschMartinBornemannUhligWall2011}. The same material is considered for the airway tree but with a stiffer Young's modulus $E=13.5$ kPa to mimic the stiffer airway walls. The air is modeled as an incompressible Newtonian fluid with density $\rho^\mathrm{f}=1$ Kg/m$^3$ and viscosity $\mu=1.5\cdot 10^{-5}$ Pa $\cdot$ s. 
The entire system is assumed to be initially at rest and the action of the thoracic muscles causing the inflation of the lungs is modeled as a displacement boundary {\color{rev1}condition}  on the outer surface of the lung lobes. This boundary data is known experimentally form a temporal sequence of CT images. As a result, the lung inflates and causes the suction of the air through the airway tree, see Figure \ref{fig:lung-sol}.

For the discretization of the problem, we consider linear tetrahedral elements for the structural domain (lung lobes and walls of the air tree) and linear tetrahedral stabilized finite elements for the air. The number of degrees of freedom for each field is given in Table \ref{table:discr:lung}. The time step is chosen as $\Delta t=1\cdot 10^{-3}$s and the problem is simulated for 750 time steps. Convergence of the Newton method is declared when the full non-linear residual fulfills $\Norm{\fmat^\mathrm{C;FSI}}_\mathrm{rms}<10^{-5}$. The stopping criterion for the GMRES iteration is $\Norm{\rmat}_2/\Norm{\rmat^0}_2<10^{-6}$. 
For the resolution of this challenging problem we consider the preconditioners SIMPLE(BGS(AMG)) and SIMPLE(AMB(BGS)) presented in Table \ref{table:preconds:lung}.

\begin{table}[ht!]
\centering
\begin{tabular}{ccccccc}
\toprule
Processors & Solid & Fluid & Grid &  Constraint     & Total & Total/Processors   \\
\midrule
64   &   1557897 &   429552 &  210498 &  5 & 2197952 & 34343\\
\bottomrule
\end{tabular}
\caption{Patient-specific lung example: Number of processors and number of degrees of for the mesh considered in this example. }
\label{table:discr:lung}
\end{table}

\begin{figure}[ht!]
\centering
\includegraphics[width=0.8\textwidth]{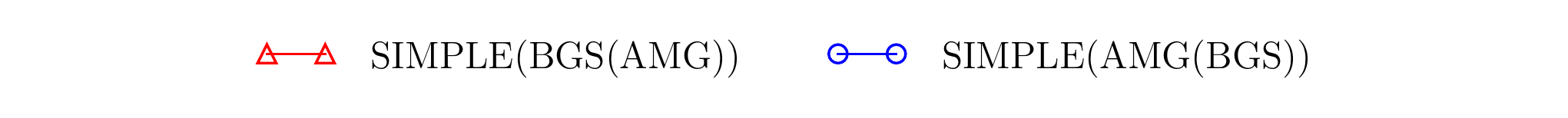}\\[0.2cm]

\begin{tabular}{cc}
\includegraphics[width=0.4\textwidth]{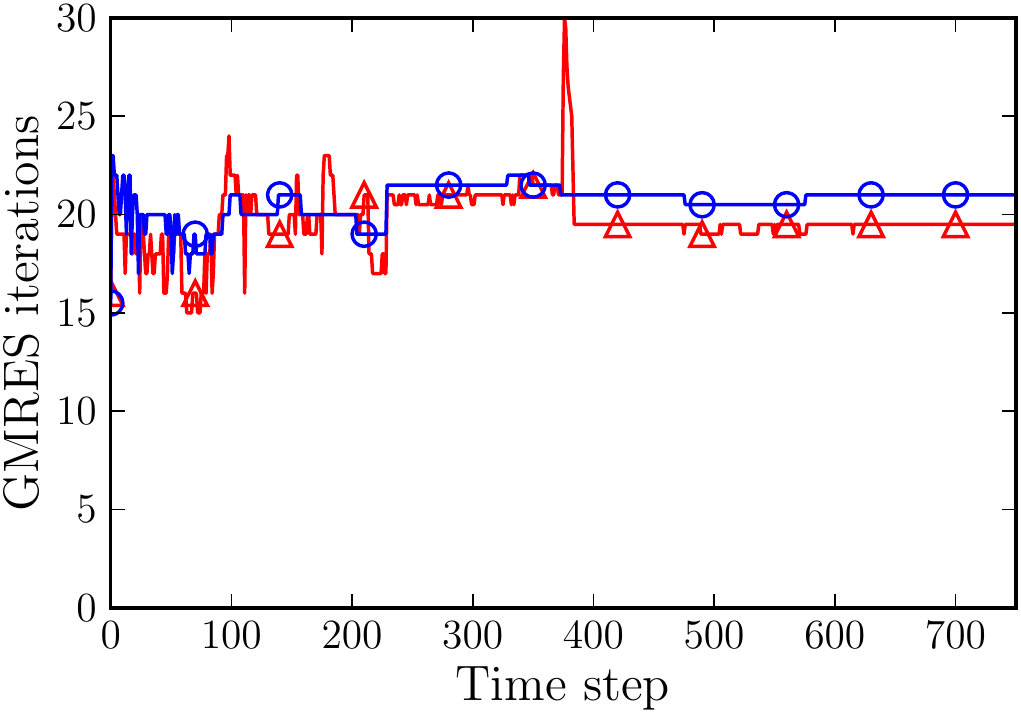} & \includegraphics[width=0.4\textwidth]{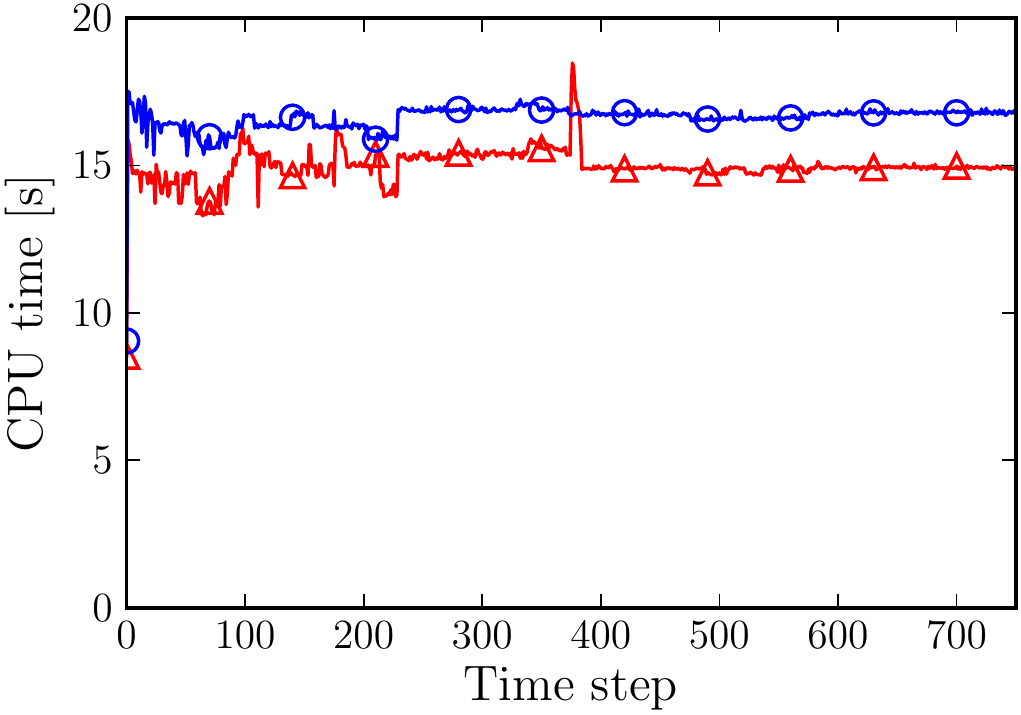}
\end{tabular}
\caption{Patient-specific lung example: Average of linear solver iterations per Newton step (left) and the corresponding CPU times (right) including the setup cost of the preconditioners.}
\label{fig:lung-iters}
\end{figure}

Figure \ref{fig:lung-iters} shows the performance of the studied  preconditioners. Note that both methods are able to solve the linear systems in about 20 iterations, which is a very good convergence behavior. This shows that the proposed methods are preconditioning the problem effectively also in this challenging {\color{rev1} real-world} example.  
The average time spent in the solution of the linear systems is  14.9 seconds for SIMPLE(BGS(AMG)) and 16.6 seconds for SIMPLE(AMG(BGS)). 
These times represent the 33.4 \% of the total simulation time for the SIMPLE(BGS(AMG)) method and 36.5 \% for the SIMPLE(AMG(BGS)). These can be considered very good results since the linear solver phase (which is usually considered the bottle neck of the computation) is taking here less than the 40\% of the total computation time.
In conclusion, the results show that the {\color{rev1} presented} methods are  effective preconditioners allowing to solve this challenging problem in decent amount of time.

\section{Concluding remarks}

{\color{rev1} Several} algebraic multigrid preconditioners have been {\color{rev1} evaluated} for solving challenging systems of linear equations resulting from the discretization of coupled problems. {\color{rev1} These methods have been considered because} they can be implemented for {\color{rev1} a} generic coupled problem and used in different applications.  The versatility of {\color{rev1} the selected techniques} is demonstrated by considering three very different coupled problems: thermo-structure interaction (TSI), fluid-structure (FSI) interaction and a complex model of the human lung. 
The numerical results show that the  proposed methods are scalable precondit{\color{rev1}i}oners in different settings including a benchmark test for TSI, a real world TSI example, and a FSI problem. Moreover, it is shown that the generic methods are efficient in terms of CPU time even though they are not optimized for a particular problem. This is showed in an FSI example by comparing the general purpose methods presented here and FSI specific preconditioners previously presented in the literature. Finally, the robustness and the wide applicability range of the methods are illustrated with a real-world patient-specific simulation of the human lung. The proposed methods lead to a very good convergence of the preconditioned GMRES method and efficient CPU times also in this challenging setting.

\section*{Acknowledgment}
{\color{rev1}
We gratefully acknowledge the support through the German Research Foundation (Deutsche Forschungsgemeinschaft -- DFG) in the framework of the Sonderforschungsbereich Transregio 40, TP D4.
} 

The authors also acknowledge the support given by the Bayerische Kompetenznetzwerk f\"ur Technisch-Wissenschaftliches Hoch- und H\"ochstleistungsrechnen (KONWIHR) in the framework of the project \emph{Efficient solvers for coupled problems in respiratory mechanics}.

The access to the high performance facilities provided by the Leibniz Supercomuting Center (Munich, Germany) in the framework of the project \emph{Coupled Problems in Computational Modeling of the Respiratory System} (project code pr32ne) is also gratefully acknowledged.

%
%
%
%
%

\bibliographystyle{unsrtnat}
\bibliography{references}

\end{document}